\documentstyle{article}

\newfont{\lat}{cmssbx9}
\newfont{\cir}{wnssdc10 scaled 900}

\def\mj{{\mathbf{1}}}
\def\pl{\!+\!}
\def\mn{\!-\!}
\def\cirk{\,{\raisebox{.3ex}{\tiny $\circ$}}\,}
\def\prop#1#2{\vspace{2ex} \noindent{\sc #1.} {\it #2} \par \vspace{2ex}}
\def\dkz{\noindent{\sc Proof. }}
\def\qed{\hfill $\dashv$}
\def\str{\rightarrow}
\def\strt{\stackrel{\textbf{.}\,}{\rightarrow}}
\def\rts{\leftarrow}

\begin{document}

\title{Coherence for Modalities}
\author{\small {\sc Kosta Do\v sen} and {\sc Zoran Petri\' c}
\\[1ex]
{\small Mathematical Institute, SANU}\\[-.5ex]
{\small Knez Mihailova 36, p.f.\ 367, 11001 Belgrade,
Serbia}\\[-.5ex]
{\small email: \{kosta, zpetric\}@mi.sanu.ac.rs}}
\date{}
\maketitle

\begin{abstract}
\noindent Positive modalities in $S4$, $S5$ and systems in their
vicinity are investigated in terms of categorial proof theory.
Coherence and maximality results are demonstrated, and connections
with mixed distributive laws and Frobenius algebras are exhibited.
\end{abstract}

\vspace{.3cm}

\noindent {\small \emph{Mathematics Subject Classification} ({\it
2010}): $\;$03B45, 03G30, 18C15, 03F07, 08C05, 18A15, 18A40,
18B10, 18C20}

\vspace{.5ex}

\noindent {\small {\it Keywords}: modality, categorial coherence,
monad, triple, comonad, adjunction, simplicial category, split
equivalence, mixed distributive laws, Frobenius monad, maximality}

\section{Introduction}

A \emph{modality} is a finite (possibly empty) sequence of the
modal operators of necessity $\Box$ and possibility $\Diamond$.
Negation is usually also allowed to occur in a modality, and the
definition we just gave would cover only \emph{positive}
modalities, but in this paper we do not consider negation (for
reasons mentioned towards the end of this introduction), and we
take \emph{modality} to be synonymous with \emph{positive
modality}. Our aim is to investigate modalities for logics in the
vicinity of $S4$ and $S5$ in terms of categorial proof theory.

The modalities in $S4$ and $S5$ are pretty well known, and one
could imagine there is nothing new to say about this topic. This
is indeed so if one wants to say just what modalities are
equivalent, and which implies which (this structure, which reduces
to just three nonequivalent positive modalities: $\Box$, the empty
modality and $\Diamond$, is very simple for $S5$). If however one
approaches this topic from the point of view of general proof
theory, or categorial proof theory, where one is interested in
identity of deductions, there are quite interesting facts about
the modalities of $S4$, $S5$ and logics in their vicinity, facts
that are not very well known, or are not known at all.

We consider deductions involving only modalities, and define
categories whose objects are these modalities, and whose arrows
may be taken as these deductions. For the logic $S4$, these happen
to be freely generated categories that have the structure of a
monad (or triple) or a comonad (for these notions, see
\cite{ML98}, Section VI.1, and Sections 3 and 8 below). The
connection between $S4$ and the notions of monad and comonad is
known (\cite{L69}, Section~1, should be the first reference for
this connection, which was exploited in particular in papers
dealing with categorial models of deductions in linear logic,
starting with \cite{SEELY89}), but here we present this matter in
a new, gradual, detailed and systematic manner, concentrating on
coherence results, some of which are presumably new. (We will, of
course, give references concerning results for which we know that
they have been previously established.)

Roughly speaking, a coherence result is a result that
characterizes a category $\cal C$, freely generated in a class of
categories, in terms of a manageable category $\cal M$. More
precisely, in a coherence result one establishes that there is a
faithful functor $G$ from $\cal C$ to $\cal M$. One may take that
$\cal C$ is syntax and $\cal M$ a model. Coherence then amounts to
proving a completeness theorem: the existence of the functor $G$
is soundness, while its faithfulness is completeness proper. As it
happens often with completeness theorems, coherence results yield
usually through the manageability of $\cal M$ an easy decision
procedure for equality of arrows in $\cal C$.

In this paper, and in general in our approach to coherence (see
\cite{DP04}), the syntactic categories $\cal C$ are indeed
constructed out of syntactic material. They correspond to logical
systems, but not to the usual systems of theorems; we have instead
systems of equations between deductions. (The usual logical
systems correspond here to the inductive definitions of terms that
stand for deductions.)

The model category $\cal M$ often has a geometrical inspiration,
and its arrows can be drawn. In this paper, the arrows of $\cal M$
will be relations of some kind, which can always be drawn. In the
first part of the paper, for categories in Sections 2-5, these
relations are either relations between finite ordinals or split
equivalences between finite ordinals. A split equivalence is an
equivalence relation on the union of two disjoint source and
target sets (see \cite{DP07}, Section 2.3, \cite{DP03a} and
\cite{DP03b}), which here we take to be finite ordinals. For the
categories in Sections 6-8, our relations are always split
equivalences between finite ordinals.

In contradistinction to coherence such as it is treated in
\cite{DP04}, the relations of $\cal M$ in this paper do not link
occurrences of propositional letters, but occurrences of modal
operators. This approach may suggest finer coherence results for
predicate logic than those obtained in \cite{DP07a}, where
quantifiers were not linked, but only predicate letters.

In this paper we deal only with modalities, which is a preparatory
matter for a treatment of wider fragments of modal logic,
involving other connectives. We believe that concerning this basic
matter we have enough new material to present, especially in
connection with $S5$, and that it is unwise to rush to wider
fragments without having settled fundamentals first. If in these
wider fragments we link both occurrences of propositional letters
and occurrences of modal operators, hoping for coherence, we enter
into a largely unchartered territory. Let us only mention that in
the presence of a lattice conjunction, which corresponds to binary
product, or a lattice disjunction, which corresponds to binary
coproduct, we should not expect straightforward coherence results
if we have both kinds of link. Problems arise with distribution of
$\Box$ over such a conjunction and distribution of $\Diamond$ over
such a disjunction. An approach to $S5$ following the present
paper would require that these distributions be isomorphisms, and
this does not square with coherence (for much the same reasons
that prevent a straightforward approach to coherence with the
isomorphic distribution of product over coproduct, which one has
in bicartesian closed categories; see \cite{DP04}, Sections 1.2
and 14.3).

The finite ordinals that are the objects of $\cal M$ in this paper
may be replaced by modalities so as to make $\cal C$ isomorphic to
a subcategory of $\cal M$ (cf.\ the parenthetical remark in the
first paragraph of Section~3 of \cite{DP08}), but when the objects
of $\cal C$ are modalities built out only of $\Box$, or only of
$\Diamond$, we need no further adjustments of $\cal M$ to have as
a consequence of coherence that $\cal C$ is isomorphic to a
subcategory of $\cal M$, and this subcategory may happen to be an
important and interesting concrete category. As an example of such
an important concrete category, we find in this paper the
simplicial category, whose arrows are the order-preserving
functions between finite ordinals, which is isomorphic to the
category $\cal C$ whose arrows may be taken as the deductions in
the modal logic $S4$ involving the modalities built out only of
$\Diamond$. As another example, we have the skeleton of the
category \emph{Finset} of finite sets, whose arrows are all the
functions between finite ordinals, which is isomorphic to the
category $\cal C$ whose arrows may be taken as the deductions
involving the modalities built out only of $\Diamond$ in an
extension of $S4$. The isomorphisms with $\cal C$ provide an
axiomatic presentation in terms of generating arrows and equations
between arrows of these important concrete categories.

Before we reach $S4$, we have in Section~2 a basic underlying
category that we call $S$. We prove for $S$ a basic simple
coherence result, which is an essential ingredient of the proofs
of coherence in the subsequent two sections (Sections 3-4) dealing
with categories related to $S4$, and in later sections. The arrows
of the category $S$ and closely related categories may be taken as
the deductions involving the modalities in the modal logics $T$
and $K4$. As a consequence of coherence for $S$, we obtain the
isomorphism of categories closely related to $S$ with respect to
the concrete categories whose arrows are respectively the
order-preserving injections and the order-preserving surjections
between finite ordinals. These isomorphisms yield axiomatic
presentations in terms of generating arrows and equations of these
concrete categories. They show also that the notions of
injectivity and surjectivity are almost the same in this context.

After $S4$ we concentrate in Section~5 on modal logics with
deductions permuting modalities. Some of these, which permute
$\Box$ with $\Box$, or $\Diamond$ with $\Diamond$, would not be
distinguished from $S4$ in ordinary modal logic, where we are
interested only in theorems, and not in deductions. From our
proof-theoretical point of view, we obtain however new logics, for
whose categories of modalities we prove coherence results. As an
interesting consequence of these results, one obtains through the
isomorphism with the category $S4_{\Diamond\chi}$ of Section~5 an
axiomatic presentation in terms of generating arrows and equations
of the skeleton of \emph{Finset} mentioned above. In this context
we also have the modal logic $S4.2$ (new from anybody's point of
view), for which we also prove coherence. The related category,
combining a monad and a comonad, was remarked independently in
attempts to describe an algebra and a coalgebra with mixed
distributive laws (see Section~5 for references).

The first part of the paper up to Section~6 is to a great extent
of an introductory character. It systematizes matters, many of
which are already known, and lays the ground for our main results
in the remainder of the paper.

In Sections 6-7 we consider categories that correspond to $S5$ and
a dual system, usually not considered in modal logic, which we
call $5S$. These categories are about combining a monad and a
comonad structure as in situations where a functor has both a left
and a right adjoint (for the notion of adjunction, see
\cite{ML98}, Section IV.1, and the beginning of Section 10 below).
These common adjoint situations do not seem to have a standard
name. In Section~8 we give them the name \emph{trijunction}, while
the corresponding monad-comonad structures, exemplified in $S5$
and $5S$, will be called \emph{dyad} and \emph{codyad}.

The dyad and codyad structures are closely related to Frobenius
algebras, a topic that has recently become rather prominent with
the proof of the equivalence between the category of commutative
Frobenius algebras and two-di\-men\-sion\-al topological quantum
field theories (see \cite{K03}). Our coherence results for the
free dyad and codyad are related to these topological results. The
difference is that with Frobenius monads, which correspond to
Frobenius algebras, $\Box$ and $\Diamond$ are not distinguished
any more, but the gist of the matter is in the results of this
paper. It is an interesting connection between modal logic and
topology, found on a different level from the well-known
connection between $S4$ operators and the topological interior and
closure operators. Here the connection with topology arises for
$S5$, and its dual $5S$. (The roots of topology and modal logic
are intermingled: one of the earliest papers in modal logic---in
some sense the first one from the modern point of view---is
\cite{K22}; in that paper Kuratowski actually introduced $S4$,
algebraically treated, for the first time, and investigated its
modalities.)

For our coherence results concerning categories that correspond to
$S5$ and $5S$, the model category $\cal M$ is a category whose
arrows are split equivalences between finite ordinals. With arrows
being relations between finite ordinals, we would obtain different
categories that correspond to $S5$ and $5S$, with which we do not
deal in this paper.

In the final sections of the paper (Sections 9-11) we deal with
the property of \emph{maximality} for our categories of
modalities. This is a kind of syntactic completeness for the
systems of equations of arrows that define these categories, a
property analogous to the Post completeness (which should be
called \emph{Bernays completeness}; see \cite{Z99}) of classical
propositional logic. Maximality is important because it shows that
not only our categories with relations, but any nontrivial
category modelling our categories of modalities could serve as a
faithful model.

Beyond our nontrivial categories, for which we have coherence and
maximality, we find preorders, where all arrows with the same
source and target, i.e.\ all deductions with the same premise and
conclusion, are equal. These categories are trivial from the point
of view of general proof theory, but it is not trivial to find
systems of equations that guarantee that they are preorders, as we
do in the sections on maximality at the end of the paper. These
are also coherence results, in the sense of the earliest coherence
result there is; namely, Mac Lane's coherence result for monoidal
categories in \cite{ML63}.

Matters pertaining to coherence for modalities involving classical
negation would not change significantly the picture we present,
and this is why we concentrate on positive modalities only. In the
presence of binary connectives, conjunction, disjunction, or
implication, where we would not deal only with modalities any
more, matters would however change considerably. The
distributivity of the necessity operator $\Box$ over conjunction,
or, dually, of the possibility operator $\Diamond$ over
disjunction, which normal modal logics require, introduces
particular problems for our model categories $\cal M$ with
relations. We leave these problems for a separate treatment.

In this paper we do not deal with categories of modalities that
correspond to Frobenius monads, where $\Box$ and $\Diamond$ are
isomorphic (they actually coincide), and where these modal
operators lose the meaning they have usually in modal logic. These
categories are very interesting, in particular because of their
relationship with topological quantum field theories mentioned
above, but we prefer not to extend further a sufficiently long
paper. We leave for \cite{DP08a} these matters, which are at the
limits of logic in the strict sense.

For the proof of our coherence results we rely on normal forms.
Although these normal forms are similar to those found in proof
theory, they are not inspired by cut elimination in the style of
Gentzen. Cut elimination however would work too, at least in some
cases (see the comments in the next section). These normal forms
may be easier to connect with natural deduction than with
Gentzen's sequent systems. The possibility to obtain these normal
forms is a proof-theoretical justification that our equations
between deductions are well chosen. Our coherence and our
maximality results provide other such justifications. (For an
extended discussion of these matters see \cite{DP04}.)

We assume for this paper an acquaintance with only rather basic
notions of category theory, which may all be found in \cite{ML98}.
Practically no knowledge of modal logic is assumed, except for the
sake of motivation, which may be gathered from \cite{HC96}. Some
further references concerning category theory and modal logic will
be given later in the paper.

\section{The category $S$}

We define in this section a basic category called $S$, and prove
for it a basic simple coherence result, which will be an essential
ingredient of the proofs of coherence in later sections. We
introduce first some terminology and notation.

Every \emph{arrow term}, i.e.\ term for an arrow in a category,
has a \emph{type} assigned to it; a type is a pair of objects
$(A,B)$ where $A$ is the \emph{source} and $B$ the \emph{target}.
We use $f,g,h,\dots,$ sometimes with indices, as variables for
arrow terms, and ${f\!:A\vdash B}$ indicates that the arrow term
$f$ is of type $(A,B)$. (The turnstile $\vdash$ reminds us here
that our arrows may be taken as deductions.)

The objects of the category $S$ are the finite ordinals. The
primitive arrow terms of $S$ are
\begin{tabbing}
\mbox{\hspace{14em}}\=$\mj_n\,$\=$: n\vdash n$,\\*[.5ex]
\>$\xi_n$\>$: n\pl 1\vdash n$.
\end{tabbing}
The arrow terms of $S$ are closed under the operations:
\begin{tabbing}
\mbox{\hspace{1.7em}}\=if ${f\!:n\vdash m}$ and ${g\!:m\vdash k}$
are arrow terms, then so is ${(g\cirk f)\!:n\vdash k}$;
\\*[1ex]
\>if ${f\!:n\vdash m}$ is an arrow term, then so is ${Mf\!:n\pl
1\vdash m\pl 1}$.
\end{tabbing}
We take for granted the outermost parentheses of arrow terms, and
omit them. (Further omissions of parentheses will be permitted by
the associativity of \cirk, namely, ({\it cat}~2) below.)

The arrows of the category $S$, and of analogous syntactic
categories considered in this paper, will be made of this
syntactic material in the manner described in detail in
\cite{DP04} (Chapter~2). The arrows of $S$ are equivalence classes
of arrow terms such that the following equations (which always
have arrow terms of the same type on the two sides of $=$) are
satisfied for ${f\!:n\vdash m}$:
\begin{tabbing}
\hspace{0em}\emph{categorial equations}:
\\*[1ex]
\hspace{1.7em}\= ({\it cat}~1)\hspace{8em}\= $f\cirk
\mj_n=\mj_m\cirk f=f$,
\\*[1ex]
\>({\it cat}~2)\>$h\cirk (g\cirk f)=(h\cirk g)\cirk f$,
\\[1.5ex]
\hspace{0em}\emph{functorial equations}:
\\*[1ex]
\>$(M1)$\>\=$M\mj_n=\mj_{n+1}$,
\\*[1ex]
\>$(M2)$\>$M(g\cirk f)=Mg\cirk Mf$,
\\[1.5ex]
\hspace{0em}\emph{naturality equation}:
\\*[1ex] \>($\xi$~{\it nat})\>\=$\xi_m\cirk Mf=f\cirk\xi_n$.
\end{tabbing}

\vspace{1ex}

The functorial equations say that $M$, where $Mn$ is ${n\pl 1}$,
is an endofunctor of $S$ (i.e.\ a functor from $S$ to $S$). The
naturality equation ($\xi$~{\it nat}) can be replaced for $S$ by
the two equations
\begin{tabbing}
\hspace{1.7em}\=({\it cat}~1)\hspace{8em}\= $f\cirk
\mj_n=\mj_m\cirk f=f$,\kill

\>($\xi^{MM}$~{\it nat})\>$\xi_{m+1}\cirk MMf=Mf\cirk\xi_{n+1}$,
\\*[1ex]
\>($\xi\; M$)\>\=$\xi_n\cirk M\xi_n=\xi_n\cirk\xi_{n+1}$.
\end{tabbing}
(For other categories, to be considered in later sections, the
last two equations will not necessarily yield ($\xi$~{\it nat}),
because of the presence of arrows different from $\xi_n$.)

The category $S$ can be presented as a strict monoidal category
(where associativity arrows are identity arrows), with tensor
product given by addition of natural numbers. So presented, it
would be a product category (PRO) without permutation in the sense
of \cite{ML65} (Chapter~V; for a more recent reference see
\cite{L04}). Many of the categories considered later in this paper
have analogously the structure of a product category, or a product
category with permutation, i.e.\ symmetry (PROP).

For ${k\geq 0}$, let $M^k$ be the sequence of $k$ occurrences of
$M$. Every arrow term of the form $M^k\xi_n$ is called a
$\xi$\emph{-term}. For ${n\geq 1}$, an arrow term of the form
${f_n\cirk\ldots\cirk f_1}$, where $f_1$ is $\mj_m$ for some $m$
and for every ${i\in\{2,\ldots,n\}}$ we have that $f_i$ is a
$\xi$-term, is called a \emph{developed} arrow term.

It is easy to show by using categorial and functorial equations
that the following lemma holds for $S$, and, with an appropriate
understanding of ``developed arrow term'', for all the categories
that will be considered in this paper.

\prop{Development Lemma}{For every arrow term $f$ there is a
developed arrow term $f'$ such that $f=f'$.}

Next we define inductively two functors, $G^\varepsilon$ and
$G^\delta$, from $S$ to the category \emph{Rel}, whose objects are
again the finite ordinals, and whose arrows are the relations
between finite ordinals; composition in \emph{Rel} is composition
of relations, and the identity arrows are identity relations. For
${\alpha\in\{\varepsilon,\delta\}}$, let $G^\alpha n$ be $n$, let
$G^\alpha\mj_n$ be the identity relation on $n$, and let
$G^\alpha\xi_0$ be the empty relation between $1$, which is equal
to $\{\emptyset\}$, and $0$, which is equal to $\emptyset$. For
${n\geq 1}$, we have clauses corresponding to the following
pictures:
\begin{center}
\begin{picture}(305,40)

\put(70,10){\circle{2}} \put(102,10){\circle{2}}
\put(54,30){\circle{2}} \put(70,30){\circle{2}}
\put(102,30){\circle{2}}

\put(80,20){$\ldots$}

\put(70,12){\line(0,1){16}} \put(102,12){\line(0,1){16}}

\put(70,7){\makebox(0,0)[t]{\scriptsize$n\mn 1$}}
\put(102.5,7){\makebox(0,0)[t]{\scriptsize$0$}}
\put(54,33){\makebox(0,0)[b]{\scriptsize$n$}}
\put(70,33){\makebox(0,0)[b]{\scriptsize$n\mn 1$}}
\put(102.5,33){\makebox(0,0)[b]{\scriptsize$0$}}
\put(35,20){\makebox(0,0)[r]{$G^\varepsilon\xi_n$}}

\put(210,10){\circle{2}} \put(242,10){\circle{2}}
\put(194,30){\circle{2}} \put(210,30){\circle{2}}
\put(242,30){\circle{2}}

\put(220,20){$\ldots$}

\put(210,12){\line(0,1){16}} \put(242,12){\line(0,1){16}}
\put(210,12){\line(-1,1){16}}

\put(210,7){\makebox(0,0)[t]{\scriptsize$n\mn 1$}}
\put(242.5,7){\makebox(0,0)[t]{\scriptsize$0$}}
\put(194,33){\makebox(0,0)[b]{\scriptsize$n$}}
\put(210,33){\makebox(0,0)[b]{\scriptsize$n\mn 1$}}
\put(242.5,33){\makebox(0,0)[b]{\scriptsize$0$}}
\put(175,20){\makebox(0,0)[r]{$G^\delta\xi_n$}}

\end{picture}
\end{center}
We have ${G^\alpha(g\cirk f)=G^\alpha g\cirk G^\alpha f}$, where
$\cirk$ on the right-hand side is composition of relations, and
for every ${f\!:n\vdash m}$ we have that the relation $G^\alpha
Mf\subseteq {(n\pl 1)\times (m\pl 1)}$ is obtained by adding the
pair $(n,m)$ to the relation $G^\alpha f\subseteq {n\times m}$.

We easily check by induction on the length of derivation that if
${f=g}$ in $S$, then ${G^\alpha f=G^\alpha g}$ in \emph{Rel}.
Hence $G^\alpha$ so defined is indeed a functor. Our purpose is to
show that the functors $G^\alpha$ are faithful functors.

A developed arrow term of $S$ is said to be in \emph{normal form}
when it has no subterm of the form ${M^k\xi_n\cirk
M^{k+l}\xi_{n-l+1}}$ for ${l\geq 1}$. That every arrow term of $S$
is equal in $S$ to an arrow term in normal form follows from the
Development Lemma and from the following equations of $S$ for
${l\geq 1}$, which for $k=0$ and $l=1$ give the equation $(\xi\;
M)$, and which could replace ($\xi$~{\it nat}) in the
axiomatization of the equations for $S$:

\begin{tabbing}
\hspace{1.7em}\=({\it cat}~1)\hspace{3.5em}\=$f\cirk
\mj_n=\mj_m\cirk f=f$,\kill

\>$(\xi\; M^l)$\>$M^k\xi_n\cirk
M^{k+l}\xi_{n-l+1}=M^{k+l-1}\xi_{n-l+1}\cirk M^k\xi_{n+1}$.
\end{tabbing}
Note that the sum of the superscripts of $M$ on the right-hand
side is strictly smaller than that sum on the left-hand side.

We can easily establish the following lemma.

\prop{Auxiliary Lemma}{If $f$ and $g$ are in normal form and
${G^\alpha f=G^\alpha g}$, then $f$ and $g$ are the same arrow
term.}

\noindent To prove this lemma we proceed by induction on the
number of $\xi$-terms in $f$ and $g$, which must be equal.

We infer immediately from the Auxiliary Lemma that the normal form
of an arrow term is unique. This fact can however easily be
established directly by confluence (i.e.\ the Church-Rosser
property) of reductions that consist in passing from the left-hand
side of $(\xi\; M^l)$ to the right-hand side.

We can then infer easily the following result.

\prop{$S$ Coherence}{The functors $G^\varepsilon$ and $G^\delta$
from $S$ to {\it Rel} are faithful.}

This coherence result could alternatively be established by
relying on a sequent presentation in the style of Gentzen (as in
\cite{G35}) of the category $S$. Instead of the primitive arrow
terms $\xi_n$ we would have the operation on arrow terms:
\begin{tabbing}
\hspace{1.7em}if ${f\!:n\vdash m}$ is an arrow term, then so is
${M_Lf\!:n\pl 1\vdash m}$,
\end{tabbing}
which is easily defined in terms of $\xi_n$, and vice versa, in
the presence of $\cirk$ and $\mj_n$; namely, we have
$M_Lf=_{df}f\cirk\xi_n$ and $\xi_n=_{df}M_L\mj_n$. The following
equations:
\[
\begin{array}{rcl}
g\cirk M_Lf & \!\!\!\!=\!\!\!\! & M_L(g\cirk f),
\\[1ex]
M_Lg\cirk Mf & \!\!\!\!=\!\!\!\! & M_L(g\cirk f),
\end{array}
\]
together with $(M2)$ and ({\it cat}~1), enable us to find for
every arrow term $f$ a composition-free arrow term $f'$ such that
${f=f'}$. The Auxiliary Lemma then holds if we replace ``in normal
form'' by ``composition-free'', and this yields $S$ Coherence.

So there are two ways to obtain a normal form. The first is to
``draw compositions out'', as we did first, and as Mac Lane does
in the Lemma of Section VII.5 of \cite{ML98} (see the next section
of this paper). The second way is to ``push compositions inside'',
until they disappear, as Gentzen would do. This is the gist of his
cut-elimination method.

Let the category $S_+$ be defined like $S$ save that we have
$\xi_n$ only for ${n\geq 1}$. It is easy to show \emph{$S_+$
Coherence}; namely, the assertion that the functors
$G^\varepsilon$ and $G^\delta$ from $S_+$ to \emph{Rel}, defined
in the same way as before, are faithful.

Let $S^{op}$ be the category opposite to $S$, and let the functor
$G^\alpha$ from $S^{op}$ to \emph{Rel} be defined by taking that
$G^\alpha f^{op}=(G^\alpha f)^{-1}$, where $R^{-1}$ is the
relation converse to $R$; on objects, $G^\alpha$ is again
identity. Then out of $S$ Coherence and $S_+$ Coherence we can
infer \emph{$S^{op}$ Coherence}, which says that these new
functors $G^\alpha$ are faithful, and \emph{$S_+^{op}$ Coherence},
which says that analogously defined functors from $S_+^{op}$ to
\emph{Rel} are faithful.

The category $S$ could be called $T_\Box$, because its arrows may
be taken as the deductions in the modal logic $T$ (the normal
modal logic with the axiom ${\Box p\str p}$ or ${p\str\Diamond
p}$, which is characterized by reflexive frames; see \cite{HC96},
p.\ 42) involving the modalities built out only of $\Box$,
provided $M$ is replaced by $\Box$. The category $S_+$ could
analogously be called $K4_\Diamond$, because its arrows may be
taken as the deductions in the modal logic $K4$ (the normal modal
logic with the axiom ${\Box p\str \Box\Box p}$ or
${\Diamond\Diamond p\str\Diamond p}$, which is characterized by
transitive frames; see \cite{HC96}, p.\ 64) involving the
modalities built out only of $\Diamond$, provided $M$ is replaced
by $\Diamond$. For analogous reasons, $S^{op}$ could be called
$T_\Diamond$, and $S_+^{op}$ could be called $K4_\Box$. The
interesting coherence results here are then the $G^\varepsilon$
instances of $S$ Coherence and $S^{op}$ Coherence, and the
$G^\delta$ instances of $S_+$ Coherence and $S_+^{op}$ Coherence,
as it will become clear in the next section.

By combining the assumptions for $T_\Box$ and $T_\Diamond$, for
two distinct modal operators $\Box$ and $\Diamond$, we would
obtain the category $T_{\Box\Diamond}$, whose arrows may be taken
as the deductions in the modal logic $T$ involving all the
positive modalities (cf.\ Section~4). We may combine analogously
the assumptions for $K4_\Box$ and $K4_\Diamond$ to obtain the
category $K4_{\Box\Diamond}$, whose arrows may be taken as the
deductions in the modal logic $K4$ involving all the positive
modalities. Since $\Box$ and $\Diamond$ do not ``cooperate'' in
$T_{\Box\Diamond}$ and $K4_{\Box\Diamond}$, we can prove easily
coherence for the first with respect to a $G^\varepsilon$ functor,
and coherence for the second with respect to a $G^\delta$ functor,
which are the interesting forms of coherence here (cf.\
Section~4).

The arrows of the category defined like $S$ save that we omit the
arrows $\xi_n$ and the equation ($\xi$~{\it nat}) may be taken as
the deductions in the minimal normal modal logic $K$ involving the
modalities built out only of $\Box$, or only of $\Diamond$. This
is however a discrete category: all its arrows are identity
arrows, and coherence for it, which is very easy to establish, is
a trivial result. (The category whose arrows may be taken as the
deductions in $K$ involving all the positive modalities is also
discrete; $\Box$ and $\Diamond$ do not cooperate in this
category.)

The $G^\varepsilon$ instance of $S^{op}$ Coherence and the
$G^\delta$ instance of $S_+$ Coherence, together with easily
established facts about the generation of order-preserving
injections and surjections between finite ordinals, yields that
$S^{op}$ is isomorphic to the category whose arrows are the
order-preserving injections between finite ordinals, and $S_+$ is
isomorphic to the category whose arrows are the order-preserving
surjections between finite ordinals. All this shows that the
notions of injectivity and surjectivity are up to duality almost
the same.

\section{The categories $S4_\Box$ and $S4_\Diamond$}

We introduce now the category $S4_\Box$, whose arrows may be taken
as the deductions in the modal logic $S4$ involving the modalities
built out only of $\Box$. We identify these modalities with their
lengths, and so we take as the objects of $S4_\Box$ not these
modalities, but the natural numbers, i.e.\ finite ordinals. The
category $S4_\Box$ is isomorphic to the category $\Delta^{op}$ for
$\Delta$ being the simplicial category, i.e.\ the category whose
arrows are the order-preserving functions between finite ordinals
(see \cite{ML98}, Section VII.5, and the end of this section). The
category $S4_\Box$ is the free comonad generated by a single
object, and the opposite category $S4_\Diamond$, isomorphic to
$\Delta$, which we will consider later in this section, is the
free monad generated by a single object (see the beginning of
Section~8).

The objects of $S4_\Box$ are the finite ordinals. The primitive
arrow terms of $S4_\Box$ are $\mj_n\!: n\vdash n$ plus
\[
\begin{array}{rcl}
\varepsilon_n^\Box\!\!\!\!\! & : & \!\!\!\! n\pl 1\vdash n,\\[1ex]
\delta_n^{\Box\Box}\!\!\!\!\! & : & \!\!\!\! n\pl 1\vdash n\pl 2.
\end{array}
\]
In the notation for comonads of \cite{ML98} (Section VI.1), our
$\varepsilon^\Box$ and $\delta^{\Box\Box}$ correspond respectively
to $\varepsilon$ and $\delta$. (We write the superscripts because
we introduce in this paper a systematic notation for comonads,
monads and their combinations; see $\varepsilon^\Diamond$ and
$\delta^{\Diamond\Diamond}$ towards the end of this section, and
also the notation of Sections 6 and 7.) In \cite{DP07a}, whose
subject matter is related to the subject matter of the present
paper, $\iota$ (derived from \emph{instantiation}) corresponds to
$\varepsilon$ as it is used in this paper.

The operations on arrow terms are as for $S$, with $M$ replaced by
$\Box$. The arrows of $S4_\Box$ are obtained by assuming the
following equations besides the categorial and functorial
equations:
\begin{tabbing}
\hspace{1.7em}\=({\it cat}~1)\hspace{8em}\= $f\cirk
\mj_n=\mj_m\cirk f=f$,\kill

\>($\varepsilon^\Box$~{\it nat})\>\=$\varepsilon_m^\Box\cirk \Box
f=f\cirk\varepsilon_n^\Box$,
\\[2ex]
\>($\delta^{\Box\Box}$~{\it nat})\>$\Box\Box
f\,$\=$\cirk\delta_n^{\Box\Box}\,$\=$=\delta_m^{\Box\Box}\cirk\Box
f$,
\\*[1ex]
\>$(\delta^{\Box\Box})$\>$\Box\delta_n^{\Box\Box}$\>$\cirk\delta_n^{\Box\Box}$\>$=
\delta_{n+1}^{\Box\Box}\cirk\delta_n^{\Box\Box}$,
\\[2ex]
\>$(\Box\Box\beta)$\>$\varepsilon_{n+1}^\Box$\=$\cirk\delta_n^{\Box\Box}\,$\=$=\mj_{n+1}$,
\\*[1ex]
\>$(\Box\Box\eta)$\>$\Box\varepsilon_n^\Box$\>$\cirk\delta_n^{\Box\Box}$\>$=\mj_{n+1}$.
\end{tabbing}

\vspace{1ex}

The naturality equation ($\varepsilon^\Box$~{\it nat}) is the
instance of ($\xi$~{\it nat}) for $\xi$ being $\varepsilon^\Box$,
while the naturality equation ($\delta^{\Box\Box}$~{\it nat}) and
the equation $(\delta^{\Box\Box})$ are obtained from the equations
($\xi^{MM}$~{\it nat}) and $(\xi\; M)$ adapted to
$\xi_{n+1}^{op}\!:n\pl 1\vdash n\pl 2$, which has the type of
$\delta_n^{\Box\Box}$. We may take that the equations for
$S4_\Box$, except the new equations $(\Box\Box\beta)$ and
$(\Box\Box\eta)$, are obtained from those for $S$ and $S_+^{op}$,
provided that in the presentation of $S_+^{op}$ we have instead of
($\xi^{op}$~{\it nat}) the two equations ($\xi^{MM}$~{\it nat})
and $(\xi\; M)$ adapted to $\xi_{n+1}^{op}$, which we mentioned
above. The equations for $S4_\Box$ above correspond exactly to the
equations for the category $\Delta^{op}$ obtained from the
equations (11), (12) and (13) for $\Delta$ in \cite{ML98} (Section
VII.5).

The functor $G$ from $S4_\Box$ to \emph{Rel} is defined by the
clauses
\[
\begin{array}{rcl}
G\varepsilon_n^\Box & \!\!\!\!\!=\!\!\!\! & G^\varepsilon\xi_n,
\\[1ex]
G\delta_n^{\Box\Box} & \!\!\!\!\!=\!\!\!\! &
(G^\delta\xi_{n+1})^{-1};
\end{array}
\]
otherwise, $G$ is defined like $G^\alpha$ from the preceding
section. These clauses correspond to the following pictures:
\begin{center}
\begin{picture}(305,40)

\put(70,10){\circle{2}} \put(102,10){\circle{2}}
\put(54,30){\circle{2}} \put(70,30){\circle{2}}
\put(102,30){\circle{2}}

\put(80,20){$\ldots$}

\put(70,12){\line(0,1){16}} \put(102,12){\line(0,1){16}}

\put(70,7.5){\makebox(0,0)[t]{\scriptsize$n\mn 1$}}
\put(102.5,7){\makebox(0,0)[t]{\scriptsize$0$}}
\put(54,33){\makebox(0,0)[b]{\scriptsize$n$}}
\put(70,33){\makebox(0,0)[b]{\scriptsize$n\mn 1$}}
\put(102.5,33){\makebox(0,0)[b]{\scriptsize$0$}}
\put(35,20){\makebox(0,0)[r]{$G\varepsilon_n^\Box$}}

\put(210,10){\circle{2}} \put(226,30){\circle{2}}
\put(226,10){\circle{2}} \put(258,10){\circle{2}}
\put(194,10){\circle{2}} \put(210,30){\circle{2}}
\put(258,30){\circle{2}}

\put(236,20){$\ldots$}

\put(210,12){\line(0,1){16}} \put(226,12){\line(0,1){16}}
\put(258,12){\line(0,1){16}} \put(210,28){\line(-1,-1){16}}

\put(210,6){\makebox(0,0)[t]{\scriptsize$n$}}
\put(226,7.5){\makebox(0,0)[t]{\scriptsize$n\mn 1$}}
\put(226,33){\makebox(0,0)[b]{\scriptsize$n\mn 1$}}
\put(258.5,7){\makebox(0,0)[t]{\scriptsize$0$}}
\put(194,7.5){\makebox(0,0)[t]{\scriptsize$n\pl 1$}}
\put(210,33){\makebox(0,0)[b]{\scriptsize$n$}}
\put(258.5,33){\makebox(0,0)[b]{\scriptsize$0$}}
\put(175,20){\makebox(0,0)[r]{$G\delta_n^{\Box\Box}$}}

\end{picture}
\end{center}
where the parts of the pictures involving ${0,\ldots,n\mn 1}$ do
not exist if ${n=0}$.

It is well known that $G$ so defined is a faithful functor (see
\cite{G58}, \emph{Appendice}, \cite{LAW69}, pp.\ 148ff,
\cite{L69}, p.\ 95, \cite{A74}, p. 10, \cite{D99}, Section 5.9,
and \cite{La03}, Section 2.2; among these references \cite{L69}
and \cite{D99} rely on Gentzen's cut-elimination method). We will
however prove this again by relying on the coherence results of
the preceding section. This proof is otherwise like Mac Lane's
proof of an analogous result in \cite{ML98} (Section VII.5).

\prop{$S4_\Box$ Coherence}{The functor $G$ from $S4_\Box$ to {\it
Rel} is faithful.}

\dkz We say that an arrow term $f_2\cirk f_1$ of $S4_\Box$ is in
\emph{normal form} when $\delta$ does not occur in $f_1$ and
$\varepsilon$ does not occur in $f_2$. By using the equations of
$S4_\Box$, it is easy to establish that every arrow term of
$S4_\Box$ is equal to an arrow term in normal form.

For $f$ and $g$ arrow terms of $S4_\Box$ of the same type, let
${Gf=Gg}$. Then ${f=f_2\cirk f_1}$ and ${g=g_2\cirk g_1}$ for
${f_2\cirk f_1}$ and ${g_2\cirk g_1}$ in normal form. So
$Gf_2\cirk Gf_1=Gg_2\cirk Gg_1$. It is easy to see that for every
arrow term $f$ of $S4_\Box$, the relation converse to $Gf$ is an
order-preserving function. Every order-preserving function
${h\!:m\str n}$ is equal to the composition ${h_2\cirk h_1\!:m\str
n}$ for a unique order-preserving surjection ${h_1\!:m\str k}$ and
a unique order-preserving injection ${h_2\!:k\str n}$, where $k$
is the cardinality of the image of $h$ (see \cite{MM92}, Section
IV.6, Propositions 1 and 2, for a more general categorial result,
with the help of which this can be inferred). For future
reference, we call this the \emph{surjection-injection}
decomposition of order-preserving functions between finite
ordinals.

We use this surjection-injection decomposition to establish that
${Gf_1=Gg_1}$ and ${Gf_2=Gg_2}$. Then we use the $G^\varepsilon$
instance of $S$ Coherence to establish that ${f_1=g_1}$, and the
$G^\delta$ instance of $S_+^{op}$ Coherence to establish that
${f_2=g_2}$, from which it follows that ${f=g}$ in $S4_\Box$. \qed

\vspace{2ex}

The normal form introduced in this proof, which is suggested by
the sur\-jec\-tion-injection decomposition, could be replaced in
our proof by a normal form suggested by another decomposition of
order-preserving functions between finite ordinals, which should
be called the \emph{injection-surjection} decomposition. In this
other decomposition we have that every order-preserving function
${h\!:m\str n}$ is equal to ${h_2\cirk h_1\!:m\str n}$ for a
unique order-preserving injection ${h_1\!:m\str k}$ and a unique
order-preserving surjection ${h_2\!:k\str n}$, where $k$ is ${m\pl
n}$ minus the cardinality of the image of $h$. This new normal
form is obtained from the previous one ${f_2\cirk f_1}$ by
applying naturality equations until we obtain ${f_1'\cirk f_2'}$
such that $\varepsilon$ does not occur in $f_2'$ and $\delta$ does
not occur in $f_1'$. The old normal form is \emph{thin}: the
cardinality of the interpolated $k$ is the least possible; the new
normal form is \emph{thick}: the cardinality of the interpolated
$k$ can now be greater than in the thin normal form, and is in a
certain sense maximal (cf.\ \cite{D99}, Section 0.3.5).

Let the category $S4_\Diamond$ be $S4_\Box^{op}$ where $\Box$ is
written $\Diamond$, while $(\varepsilon_n^\Box)^{op}\!:n\vdash
n\pl 1$ and $(\delta_n^{\Box\Box})^{op}\!:n\pl 2\vdash n\pl 1$ are
written
\[
\begin{array}{rcl}
\varepsilon_n^\Diamond\!\!\!\!\! & : & \!\!\!\! n\vdash n\pl 1,\\[1ex]
\delta_n^{\Diamond\Diamond}\!\!\!\!\! & : & \!\!\!\! n\pl 2\vdash
n\pl 1,
\end{array}
\]
respectively. (In the notation for monads of \cite{ML98}, Section
VI.1, our $\varepsilon^\Diamond$ and $\delta^{\Diamond\Diamond}$
correspond respectively to $\eta$ and $\mu$.) The arrows of the
category $S4_\Diamond$ may be taken as the deductions in the modal
logic $S4$ involving the modalities built out only of $\Diamond$.

Let the functor $G$ from $S4_\Diamond$ to \emph{Rel} be defined by
taking that $Gf^{op}=_{df}(Gf)^{-1}$, where on the right-hand side
$G$ is the functor from $S4_\Box$ to \emph{Rel}; on objects, $G$
is identity. This means that we have clauses corresponding to the
following pictures, obtained from the pictures given above for
$G\varepsilon_n^\Box$ and $G\delta_n^{\Box\Box}$ by putting them
upside down (and taking for granted the line involving ${n\mn 1}$
in the right picture):
\begin{center}
\begin{picture}(305,40)

\put(70,10){\circle{2}} \put(102,10){\circle{2}}
\put(54,10){\circle{2}} \put(70,30){\circle{2}}
\put(102,30){\circle{2}}

\put(80,20){$\ldots$}

\put(70,12){\line(0,1){16}} \put(102,12){\line(0,1){16}}

\put(70,7){\makebox(0,0)[t]{\scriptsize$n\mn 1$}}
\put(102.5,7){\makebox(0,0)[t]{\scriptsize$0$}}
\put(54,5.5){\makebox(0,0)[t]{\scriptsize$n$}}
\put(70,33){\makebox(0,0)[b]{\scriptsize$n\mn 1$}}
\put(102.5,33){\makebox(0,0)[b]{\scriptsize$0$}}
\put(35,20){\makebox(0,0)[r]{$G\varepsilon_n^\Diamond$}}

\put(210,10){\circle{2}} \put(242,10){\circle{2}}
\put(194,30){\circle{2}} \put(210,30){\circle{2}}
\put(242,30){\circle{2}}

\put(220,20){$\ldots$}

\put(210,12){\line(0,1){16}} \put(242,12){\line(0,1){16}}
\put(210,12){\line(-1,1){16}}

\put(210,7){\makebox(0,0)[t]{\scriptsize$n$}}
\put(242.5,7){\makebox(0,0)[t]{\scriptsize$0$}}
\put(194,33){\makebox(0,0)[b]{\scriptsize$n\pl 1$}}
\put(210,33){\makebox(0,0)[b]{\scriptsize$n$}}
\put(242.5,33){\makebox(0,0)[b]{\scriptsize$0$}}
\put(175,20){\makebox(0,0)[r]{$G\delta_n^{\Diamond\Diamond}$}}

\end{picture}
\end{center}
Then out of $S4_\Box$ Coherence we can infer $S4_\Diamond$
\emph{Coherence}, which says that this new functor $G$ is
faithful.

This faithfulness result, together with the surjection-injection
decomposition of order-preserving functions between finite
ordinals, and the isomorphisms involving $S^{op}$ and $S_+$
mentioned at the end of the preceding section, yield that
$S4_\Diamond$ is isomorphic to the category whose arrows are the
order-preserving functions between finite ordinals, i.e.\ the
simplicial category $\Delta$.

\section{The category $S4_{\Box\Diamond}$}

We introduce now the category $S4_{\Box\Diamond}$, whose arrows
may be taken as the deductions in the modal logic $S4$ involving
all the positive modalities; namely, all the modalities built out
of both $\Box$ and $\Diamond$. The category $S4_{\Box\Diamond}$
will have the structures of a comonad and a monad.

The objects of $S4_{\Box\Diamond}$ are finite (possibly empty)
sequences of $\Box$ and $\Diamond$, sequences that we call
\emph{modalities}, and denote by $A,B,C,\ldots$ The primitive
arrow terms of $S4_{\Box\Diamond}$ are
\begin{tabbing}
\hspace{14em}$\mj_A\!: A\vdash A$,\\[1.5ex]
\hspace{8em}\=$\varepsilon_A^\Box\;\;\,$\=$: \Box A\vdash
A$,\hspace{5em}\=$\varepsilon_A^\Diamond\;\;\,$\=$:A\vdash\Diamond
A$,
\\[1ex]
\>$\delta_A^{\Box\Box}$\>$: \Box A\vdash \Box\Box
A$,\>$\delta_A^{\Diamond\Diamond}$\>$:\Diamond\Diamond
A\vdash\Diamond A$.
\end{tabbing}

The operations on the arrow terms of $S4_{\Box\Diamond}$ are
defined like the operations on the arrow terms of the category $S$
in Section~2, save that $n$, $m$ and $k$ are replaced respectively
by $A$, $B$ and $C$, while ${n\pl 1}$ and ${m\pl 1}$ are replaced
respectively by $MA$ and $MB$, where $M$ stands, as in the
preceding section, either for $\Box$ or for~$\Diamond$.

The arrows of $S4_{\Box\Diamond}$ satisfy the categorial and
functorial equations of Section~2, provided we make the
replacements just mentioned. We have moreover the equations taken
over from $S4_\Box$ and $S4_\Diamond$; namely, the equations
($\varepsilon^\Box$~{\it nat}), ($\delta^{\Box\Box}$~{\it nat}),
$(\delta^{\Box\Box})$, $(\Box\Box\beta)$ and $(\Box\Box\eta)$, and
the equations for $S4_\Diamond$ dual to these where $\Box$ is
replaced by $\Diamond$. (Some of these equations of $S4_\Diamond$
are mentioned in Section~6 when we give the equations for
$S5_{\Box\Diamond}$.) This concludes the definition of the
equations for $S4_{\Box\Diamond}$. Note that in these equations
$\Box$ and $\Diamond$ do not ``cooperate''.

We define a functor $G$ from $S4_{\Box\Diamond}$ to \emph{Rel} by
stipulating first that $GA$ is the length of the object $A$. For
$\alpha\in
I=\{\varepsilon^\Box,\delta^{\Box\Box},\varepsilon^\Diamond,\delta^{\Diamond\Diamond}\}$
and $GA=n$, let $G\alpha_A$ be defined like $G\alpha_n$, where $G$
in $G\alpha_n$ is either $G$ from $S4_\Box$ to \emph{Rel} or $G$
from $S4_\Diamond$ to \emph{Rel} (see the preceding section);
otherwise, $G$ is defined like $G^\alpha$ in Section~2. We are now
going to prove the following.

\prop{$S4_{\Box\Diamond}$ Coherence}{The functor $G$ from
$S4_{\Box\Diamond}$ to {\it Rel} is faithful.}

\dkz We say that an arrow term $f$ of $S4_{\Box\Diamond}$ is an
$\alpha$ \emph{arrow term} when no $\beta\in I-\{\alpha\}$ occurs
in $f$. The equations of $S4_{\Box\Diamond}$ enable us to find for
every arrow term $f$ of $S4_{\Box\Diamond}$ an arrow term equal to
$f$ in the \emph{normal form} $f_4\cirk f_3\cirk f_2\cirk f_1$
where $f_1$ is an $\varepsilon^\Box$ arrow term, $f_2$ is a
$\delta^{\Diamond\Diamond}$ arrow term, $f_3$ is a
$\delta^{\Box\Box}$ arrow term, and $f_4$ is an
$\varepsilon^\Diamond$ arrow term.

Suppose now that for $f$ and $g$ arrow terms of
$S4_{\Box\Diamond}$ of the same type we have ${Gf=Gg}$. For
$f_4\cirk f_3\cirk f_2\cirk f_1$ and $g_4\cirk g_3\cirk g_2\cirk
g_1$ being respectively the normal forms of $f$ and $g$, it is
easy to see that $f_i$ and $g_i$, for $1\leq i \leq 4$, are of the
same type, and that ${Gf_i=Gg_i}$.

Roughly speaking, $Gf_4$ and $Gg_4$ tell us which occurrences of
$\Diamond$ in the target of $f_4$ and $g_4$ \emph{disappear} in
their sources, and since according to $Gf=Gg$ the same of these
occurrences disappear, these sources must be the same, as well as
$Gf_4$ and $Gg_4$. So the targets of $f_3$ and $g_3$ are the same.
Since $G(f_4\cirk f_3\cirk f_2\cirk f_1)=G(g_4\cirk g_3\cirk
g_2\cirk g_1)$ and $Gf_4=Gg_4$, while $Gf_4$ and $Gg_4$ are
one-one functions, and hence left cancellable (see \cite{ML98},
Section I.5), we conclude that $G(f_3\cirk f_2\cirk
f_1)=G(g_3\cirk g_2\cirk g_1)$.

Roughly speaking, $Gf_3$ and $Gg_3$ tell us which occurrences of
$\Box$ in the target of $f_3$ and $g_3$ are \emph{amalgamated} in
their sources, and since according to $G(f_3\cirk f_2\cirk
f_1)=G(g_3\cirk g_2\cirk g_1)$ the same of these occurrences are
amalgamated, these sources must be the same, as well as $Gf_3$ and
$Gg_3$. We reason analogously in the two remaining dual cases,
where $i$ is $1$ and $2$, starting from the source of $f$ and $g$.

Then we can conclude out of $S$ and $S_+$ Coherence, and their
\emph{op} variants, that ${f_i=g_i}$, from which it follows that
${f=g}$ in $S4_{\Box\Diamond}$. \qed

\vspace{2ex}

Note that the normal form in this proof is \emph{thin}, in the
sense that the target of $f_2$, which is also the source of $f_3$,
is a minimal interpolant for decomposing $f$. Various other
thicker normal forms, with interpolants being modalities of
greater length, can be envisaged (among these there is a thickest
one). A thicker normal form, for which we will find analogues
later (see the normal forms for $S4_{\Box\Diamond\chi}$ and
$S4.2_{\Box\Diamond}$ in the next section), is $f_4\cirk f_2\cirk
f_3\cirk f_1$ for $f_i$ being as above. Note that $\Diamond$ does
not occur in the superscripts of $f_3\cirk f_1$, which hence
becomes an arrow term of $S4_\Box$ when $\Diamond$ is replaced by
$\Box$, while $\Box$ does not occur in the superscripts of
$f_4\cirk f_2$, which hence becomes an arrow term of $S4_\Diamond$
when $\Box$ is replaced by $\Diamond$.

In every situation where we have an \emph{endoadjunction}, i.e.,
where we have two functors $F$ and $G$ from a category $\cal A$ to
$\cal A$ such that $F$ is left adjoint to $G$ (for the notion of
adjunction, see the beginning of Section 10), the composite
functors $FG$ and $GF$, for $FG$ being $\Box$ and $GF$ being
$\Diamond$, together with the associated natural transformations
$\varepsilon^M$ and $\delta^{MM}$, defined in terms of the
canonical arrows of the adjunction (as in \cite{ML98}, Section
VI.1), have the structure of $S4_{\Box\Diamond}$.

\section{Modalities and permutation}

We will envisage in this section categories with arrows that
permute modalities, whose image by the functor $G$ will correspond
to the picture
\begin{center}
\begin{picture}(20,40)

\put(0,10){\line(1,1){20}} \put(20,10){\line(-1,1){20}}

\put(0,7){\makebox(0,0)[t]{$M_2$}}
\put(23,7){\makebox(0,0)[t]{$M_1$}}
\put(0,33){\makebox(0,0)[b]{$M_1$}}
\put(23,33){\makebox(0,0)[b]{$M_2$}}
\end{picture}
\end{center}
Here, $M_1$ and $M_2$, which are either $\Box$ or $\Diamond$, may
be either equal or not.

The category $S_\chi$ is defined like $S$ of Section~2, where $M$
is $\Box$ and $\xi$ is $\varepsilon^\Box$, with the additional
primitive arrow terms
\[
\chi_n^{\Box\Box}\!:\Box\Box n\vdash \Box\Box n,
\]
with $n\geq 0$, for which we assume the additional equations
\begin{tabbing}
\hspace{1.7em}\=({\it cat}~1)\hspace{7.4em}\= $f\cirk
\mj_n=\mj_m\cirk f=f$,\kill

\>($\chi^{\Box\Box}$~{\it nat})\>$\Box\Box
f\cirk\chi_n^{\Box\Box}=\chi_m^{\Box\Box}\cirk \Box\Box f$,
\\[1ex]
\>$(\chi\chi\Box)$\>$\chi_n^{\Box\Box}\cirk
\chi_n^{\Box\Box}=\mj_{n+2}$,
\\[1ex]
\>$(\chi\chi\chi\Box)$\>\=$\chi_{n+1}^{\Box\Box}\cirk
\Box\chi_n^{\Box\Box}\cirk\chi_{n+1}^{\Box\Box}
=\Box\chi_n^{\Box\Box}\cirk\chi_{n+1}^{\Box\Box}\cirk\Box\chi_n^{\Box\Box}$,
\\[2.5ex]
\>$(\varepsilon^\Box\chi^{\Box\Box})$\>\=$\varepsilon_{n+1}^\Box\cirk
\chi_n^{\Box\Box}=\Box\varepsilon_n^\Box$.
\end{tabbing}
The first three of these four equations are analogous to the
equations commonly used to present symmetric groups (see
\cite{CM57}, Section 6.2).

We define the functor $G$ from $S_\chi$ to \emph{Rel} like
$G^\varepsilon$ in Section~2 with an additional clause for
$\chi_n^{\Box\Box}$ that corresponds to the following picture:
\begin{center}
\begin{picture}(180,40)

\put(70,10){\circle{2}} \put(86,10){\circle{2}}
\put(102,10){\circle{2}} \put(136,10){\circle{2}}
\put(70,30){\circle{2}} \put(86,30){\circle{2}}
\put(102,30){\circle{2}} \put(136,30){\circle{2}}

\put(113,20){$\ldots$}

\put(70,12){\line(1,1){16}} \put(86,12){\line(-1,1){16}}
\put(102,12){\line(0,1){16}} \put(136,12){\line(0,1){16}}

\put(70,7){\makebox(0,0)[t]{\scriptsize$n\pl 1$}}
\put(86,5.5){\makebox(0,0)[t]{\scriptsize$n$}}
\put(102,7){\makebox(0,0)[t]{\scriptsize$n\mn 1$}}
\put(136.5,7){\makebox(0,0)[t]{\scriptsize$0$}}
\put(70,33){\makebox(0,0)[b]{\scriptsize$n\pl 1$}}
\put(86,33){\makebox(0,0)[b]{\scriptsize$n$}}
\put(102,33){\makebox(0,0)[b]{\scriptsize$n\mn 1$}}
\put(136.5,33){\makebox(0,0)[b]{\scriptsize$0$}}

\put(50,20){\makebox(0,0)[r]{$G\chi_n$}}

\end{picture}
\end{center}
where the part of the picture involving ${0,\ldots,n\mn 1}$ does
not exist if ${n=0}$.

To show that this functor $G$ is faithful, i.e.\ to show $S_\chi$
\emph{Coherence}, we establish first that every arrow term $f$ of
$S_\chi$ is equal in $S_\chi$ to an arrow term in the normal form
${f_2\cirk f_1}$ where $\chi^{\Box\Box}$ does not occur in $f_1$
and $\varepsilon^\Box$ does not occur in $f_2$. Here $f_1$ is an
arrow term of $S$ (with $M$ being $\Box$ and $\xi$ being
$\varepsilon^\Box$), while $f_2$ should be called a
$\chi^{\Box\Box}$ \emph{arrow term}. Note that $Gf$ determines
uniquely $Gf_1$ and $Gf_2$, as well as the target of $f_1$, which
is also the source of $f_2$. To obtain $S_\chi$ Coherence we rely
then on $S$ Coherence and on the standard presentation of
symmetric groups mentioned above, which we call \emph{Symmetric
Coherence}.

The category $S_{+\chi}^{op}$ is defined like $S_+^{op}$ of
Section~2, where $M$ is $\Box$ and $\xi^{op}$ is
$\delta^{\Box\Box}$, with ($\delta^{\Box\Box}$~{\it nat}) and
$(\delta^{\Box\Box})$ assumed instead of ($\xi^{op}$~{\it nat});
we have the additional primitive arrow terms $\chi_n^{\Box\Box}$
for which we assume the additional equations
($\chi^{\Box\Box}$~{\it nat}), $(\chi\chi\Box)$,
$(\chi\chi\chi\Box)$ and
\begin{tabbing}
\hspace{1.7em}\=({\it cat}~1)\hspace{6.5em}\= $f\cirk
\mj_n=\mj_m\cirk f=f$,\kill

\>$(\delta^{\Box\Box}\chi^{\Box\Box})$\>$\delta_{n+1}^{\Box\Box}\cirk
\chi_n^{\Box\Box}$\=$=\Box\chi_n^{\Box\Box}\cirk
\chi_{n+1}^{\Box\Box}\cirk\Box\delta_n^{\Box\Box}$,
\\[1ex]
\>$(\chi^{\Box\Box}\delta^{\Box\Box})$\>$\chi_n^{\Box\Box}\cirk
\delta_n^{\Box\Box}$\>$=\delta_n^{\Box\Box}$
\end{tabbing}
(these equations, as well as $(\varepsilon^\Box\chi^{\Box\Box})$
above, may be found in \cite{Bu93}, Section 2.2, \cite{La95},
\cite{Ma97}, \cite{G01}, Section~2, \cite{La03}, Section 2.3, and
\cite{K03}, p.\ 194).

We define the functor $G$ from $S_{+\chi}^{op}$ to \emph{Rel} as
$G^\delta$ in Section~2 with an additional clause for
$\chi_n^{\Box\Box}$ as above. To show that this functor is
faithful, i.e.\ to show $S_{+\chi}^{op}$ \emph{Coherence}, we
establish first that every arrow term $f$ of $S_{+\chi}^{op}$ is
equal in $S_{+\chi}^{op}$ to an arrow term in the normal form
${f_2\cirk f_1}$ where $\chi^{\Box\Box}$ does not occur in $f_1$
and $\delta^{\Box\Box}$ does not occur in $f_2$. So $f_1$ is an
arrow term of $S_+^{op}$, and $f_2$ is a $\chi^{\Box\Box}$ arrow
term. Note that $Gf$ determines uniquely $Gf_1$ and the target of
$f_1$. On the other hand, $Gf_2$ is not determined uniquely by
$Gf$. There is however a unique $Gf_2$ such that the number of
inversions in the permutation $Gf_2$ is minimal. We can omit all
the inversions involving $i$ and $j$ such that applying the
function $(Gf_1)^{-1}$ to $i$ and $j$ gives the same value. By
relying on that, together with the equation
$(\chi^{\Box\Box}\delta^{\Box\Box})$, we can finish our proof of
$S_{+\chi}^{op}$ Coherence by appealing to $S_+^{op}$ Coherence
and Symmetric Coherence.

The category $S_\chi^{op}$ is isomorphic to the category whose
arrows are arbitrary injections between finite ordinals, while
$S_{+\chi}$, which is $(S_{+\chi}^{op})^{op}$, is isomorphic to
the category whose arrows are arbitrary surjections between finite
ordinals. For that we rely on the possibility to decompose every
such injection into a bijection followed by an order-preserving
injection, and the same when ``injection'' is replaced by
``surjection''.

The category $S4_{\Box\chi}$ is obtained by combining what we have
assumed for $S_\chi$ and $S_{+\chi}^{op}$, as $S4$ may be obtained
from $S$ and $S_+^{op}$, and the functor $G$ from $S4_{\Box\chi}$
to \emph{Rel} is obtained by combining what we have assumed for
the functors $G$ from $S_\chi$ and $S_{+\chi}^{op}$. This category
is interesting because, as we will see below, its opposite
category is isomorphic to the category whose arrows are arbitrary
functions between finite ordinals.

We can prove $S4_{\Box\chi}$ \emph{Coherence} with respect to the
functor $G$ we have just defined by relying for every arrow term
$f$ of $S4_{\Box\chi}$ on its normal form ${f_2\cirk f_1}$ where
$f_1$ is an arrow term of $S$ and $f_2$ is an arrow term of
$S_{+\chi}^{op}$. Then we apply $S$ Coherence and $S_{+\chi}^{op}$
Coherence (for related proofs of $S4_{\Box\chi}$ Coherence see
\cite{Bu93}, Section 2.2, \cite{La95}, \cite{G01}, \cite{La03},
Section 2.3, and \cite{K03}, p.\ 195).

The category $S4_{\Diamond\chi}$ is $S4_{\Box\chi}^{op}$. We use
for $S4_{\Diamond\chi}$ the same notation that we used for
$S4_\Diamond$ (see Section~3), and we write
$\chi_n^{\Diamond\Diamond}$ for $(\chi_n^{\Box\Box})^{op}$. From
$S4_{\Box\chi}$ Coherence we can infer, of course,
$S4_{\Diamond\chi}$ \emph{Coherence}, for an appropriately defined
functor $G$ whose definition extends the definition of the functor
$G$ from $S4_\Diamond$ to \emph{Rel} (see Section~3) with a clause
for $\chi_n^{\Diamond\Diamond}$ like the clause for
$\chi_n^{\Box\Box}$.

As $S4_\Diamond$ is isomorphic to the category whose arrows are
the order-preserving functions between finite ordinals, so
$S4_{\Diamond\chi}$ is isomorphic to the category whose arrows are
arbitrary functions between finite ordinals, which is, up to
isomorphism, the skeleton of the category \emph{Finset} of finite
sets. This is shown by relying on the decomposition of every such
function into a bijection followed by an order-preserving
function.

We define the category $S4_{\Box\Diamond\chi}$ like
$S4_{\Box\Diamond}$ with the additional arrows $\chi_A^{\Box\Box}$
and $\chi_A^{\Diamond\Diamond}$, for which we have the equations
we have assumed for $S4_{\Box\chi}$ and $S4_{\Diamond\chi}$ with
$n$ and $m$ replaced respectively by $A$ and $B$, while ${n\pl 1}$
and ${n\pl 2}$ are replaced respectively by $\Box A$ and $\Box\Box
A$, or $\Diamond A$ and $\Diamond\Diamond A$, as appropriate. The
definition of the functor $G$ from $S4_{\Box\Diamond\chi}$ to
\emph{Rel} extends the definition of the functor $G$ from
$S4_{\Box\Diamond}$ to \emph{Rel} (see the preceding section) with
the clauses for $\chi_A^{\Box\Box}$ and
$\chi_A^{\Diamond\Diamond}$. Since $\Box$ and $\Diamond$ do not
``cooperate'' in $S4_{\Box\Diamond\chi}$, we can rely on
$S4_{\Box\chi}$ Coherence and $S4_{\Diamond\chi}$ Coherence to
establish $S4_{\Box\Diamond\chi}$ \emph{Coherence} with respect to
this functor~$G$. The normal form on which we rely here is
analogous to the thicker normal form for $S4_{\Box\Diamond}$
mentioned in the penultimate paragraph of the preceding section.

In the category $S4.2_{\Box\Diamond}$, which we are now going to
define, $\Box$ and $\Diamond$ will ``cooperate'' for the first
time. This category is obtained by extending what we have assumed
for $S4_{\Box\Diamond}$ with the additional primitive arrow terms
\[
\chi_A^{\Diamond\Box}\!:\Diamond\Box A\vdash\Box\Diamond A
\]
for which we assume the additional equations
\begin{tabbing}
\hspace{1.7em}\=({\it cat}~1)\hspace{7.5em}\= $f\cirk
\mj_n=\mj_m\cirk f=f$,\kill

\>($\chi^{\Diamond\Box}$~{\it nat})\>$\Box\Diamond
f\cirk\chi_A^{\Diamond\Box}=\chi_B^{\Diamond\Box}\cirk
\Diamond\Box f$,
\\[2ex]
\>$(\varepsilon^\Box\chi^{\Diamond\Box})$\>$\varepsilon_{\Diamond
A}^\Box\,$\=$\cirk
\chi_A^{\Diamond\Box}$\=$=\Diamond\varepsilon_A^\Box$,
\\*[1ex]
\>$(\varepsilon^\Diamond\chi^{\Diamond\Box})$\>$\chi_A^{\Diamond\Box}$\>$\cirk
\varepsilon_{\Box A}^\Diamond$\>$=\Box\varepsilon_A^\Diamond$,
\\[2ex]
\>$(\delta^{\Box\Box}\chi^{\Diamond\Box})$\>$\delta_{\Diamond
A}^{\Box\Box}$\>$\cirk \chi_A^{\Diamond\Box}$\=$=
\Box\chi_A^{\Diamond\Box}$\=$\cirk\chi_{\Box
A}^{\Diamond\Box}$\=$\cirk \Diamond\delta_A^{\Box\Box}$,
\\*[1ex]
\>$(\delta^{\Diamond\Diamond}\chi^{\Diamond\Box})$\>
$\chi_A^{\Diamond\Box}$\>$\cirk\delta_{\Box
A}^{\Diamond\Diamond}$\>$=
\Box\delta_A^{\Diamond\Diamond}$\>$\cirk\chi_{\Diamond
A}^{\Diamond\Box}$\>$\cirk\Diamond\chi_A^{\Diamond\Box}$,
\end{tabbing}
analogous to the equations ($\chi^{\Box\Box}$~{\it nat}),
$(\varepsilon^\Box\chi^{\Box\Box})$ and
$(\delta^{\Box\Box}\chi^{\Box\Box})$ above. The arrows of
$S4.2_{\Box\Diamond}$ may be taken as the deductions involving the
positive modalities in the modal system $S4.2$, which extends $S4$
with a principle corresponding to the type of
$\chi_A^{\Diamond\Box}$ (see \cite{HC96}, p.\ 134). All the
equations assumed above for $S4.2_{\Box\Diamond}$ may be found in
\cite{W08} (Section 5.3), in connection with \emph{mixed
distributive} or \emph{entwining} natural transformations in
structures that combine a comonad and a monad; these distributive
laws stem from \cite{B69}.

We define the functor $G$ from $S4.2_{\Box\Diamond}$ to \emph{Rel}
as the functor $G$ from $S4_{\Box\Diamond}$ to \emph{Rel} with an
additional clause for $\chi_A^{\Diamond\Box}$ that corresponds to
the picture given above for $G\chi_n^{\Box\Box}$. We can then show
the following.

\prop{$S4.2_{\Box\Diamond}$ Coherence}{The functor $G$ from
$S4.2_{\Box\Diamond}$ to {\it Rel} is faithful.}

\dkz We establish first that every arrow term $f$ of
$S4.2_{\Box\Diamond}$ is equal in $S4.2_{\Box\Diamond}$ to an
arrow term in the \emph{normal form} ${f_3\cirk f_2\cirk f_1}$
where $\varepsilon^\Diamond$, $\delta^{\Diamond\Diamond}$ and
$\chi^{\Diamond\Box}$ do not occur in $f_1$, while
$\varepsilon^\Box$, $\delta^{\Box\Box}$ and $\chi^{\Diamond\Box}$
do not occur in $f_3$, and $\varepsilon^\Box$,
$\delta^{\Box\Box}$, $\varepsilon^\Diamond$ and
$\delta^{\Diamond\Diamond}$ do not occur in $f_2$. If we replace
$\Diamond$ by $\Box$, then $f_1$, in whose superscripts $\Diamond$
does not occur, becomes an arrow term of $S4_\Box$, and if we
replace $\Box$ by $\Diamond$, then $f_3$, in whose superscripts
$\Box$ does not occur, becomes an arrow term of $S4_\Diamond$.
This normal form is analogous to the thicker normal form for
$S4_{\Box\Diamond}$ mentioned in the penultimate paragraph of the
preceding section. Note that $Gf$ determines uniquely $Gf_1$,
$Gf_2$ and $Gf_3$, as well as the targets of $f_1$ and $f_2$. Then
we can apply $S4_\Box$ Coherence and $S4_\Diamond$ Coherence as
far as $f_1$ and $f_3$ are concerned. As far as $f_2$ is
concerned, we can establish an easy coherence result for
categories that involve only the $\chi_A^{\Diamond\Box}$ arrows
and the functors $\Box$ and $\Diamond$, and where we have only the
categorial and functorial equations and the naturality equation
($\chi^{\Diamond\Box}$~{\it nat}). (This is accomplished by a
confluence technique; cf.\ \cite{DP04}.)\qed

\vspace{2ex}

We could define a category analogous to $S4.2_{\Box\Diamond}$ that
would have instead of the arrows $\chi_A^{\Diamond\Box}$ the
arrows with converse types:
\[
\chi_A^{\Box\Diamond}\!:\Box\Diamond A\vdash \Diamond\Box A,
\]
and appropriate equations analogous to those of
$S4.2_{\Box\Diamond}$, which deliver coherence. The equations
involving explicitly $\chi^{\Box\Diamond}$ are obtained from the
equations ($\chi^{\Diamond\Box}$~{\it nat}),
$(\varepsilon^\Box\chi^{\Diamond\Box})$,
$(\varepsilon^\Diamond\chi^{\Diamond\Box})$,
$(\delta^{\Box\Box}\chi^{\Diamond\Box})$ and
$(\delta^{\Diamond\Diamond}\chi^{\Diamond\Box})$ by taking
$\chi_A^{\Box\Diamond}$ as the inverse of $\chi_A^{\Diamond\Box}$.
For example, from $(\delta^{\Box\Box}\chi^{\Diamond\Box})$ we
obtain the equation
\[
\chi_{\Box A}^{\Box\Diamond}\cirk\Box\chi_A^{\Box\Diamond}\cirk
\delta_{\Diamond A}^{\Box\Box}=
\Diamond\delta_A^{\Box\Box}\cirk\chi_A^{\Box\Diamond}.
\]

The arrows of this category may be taken as the deductions
involving the positive modalities in McKinsey's modal system
$S4.1$, also called $S4M$ (for historical comments see
\cite{HC96}, p.\ 143, note 7), whose theorems are not included in
$S5$. Coherence for this category is demonstrated quite
analogously to what we had for $S4.2_{\Box\Diamond}$.

We can also envisage the category with both
$\chi_A^{\Diamond\Box}$ and $\chi_A^{\Box\Diamond}$ arrows, which
would be isomorphisms inverse to each other. Coherence for that
category is again shown analogously. To this last category we can
also add the arrows $\chi_A^{\Box\Box}$ and
$\chi_A^{\Diamond\Diamond}$, and again obtain easily a coherence
result.

\section{The category $S5_{\Box\Diamond}$}

We introduce now the category $S5_{\Box\Diamond}$, whose arrows
may be taken as the deductions in the modal logic $S5$ involving
the positive modalities. As $S4_{\Box\Diamond}$, this category
will have the structures of a comonad and a monad, which however
will now ``cooperate''.

We define the category $S5_{\Box\Diamond}$ like the category
$S4_{\Box\Diamond}$ with the following additions. We have the
additional primitive arrow terms
\[
\delta_A^{\Box\Diamond}\!:\Diamond A\vdash \Box\Diamond A,
\hspace{5em} \delta_A^{\Diamond\Box}\!:\Diamond\Box A\vdash \Box
A.
\]
We use $\delta_A^{\Box M}$ for either $\delta_A^{\Box\Box}$ or
$\delta_A^{\Box\Diamond}$, and likewise $\delta_A^{\Diamond M}$
for either $\delta_A^{\Diamond\Diamond}$ or
$\delta_A^{\Diamond\Box}$. The equations of $S5_{\Box\Diamond}$
are obtained by assuming those assumed for $S4_{\Box\Diamond}$ and
the following additional equations:
\begin{tabbing}
\hspace{.4em}\=($\delta^{\Box M}$~{\it nat})\hspace{.5em}$\Box
Mf\cirk\delta^{\Box M}_A=\delta^{\Box M}_B\cirk Mf$,\hspace{1.5em}
($\delta^{\Diamond M}$~{\it nat})\hspace{.5em} $\delta^{\Diamond
M}_B\cirk \Diamond Mf=Mf\cirk\delta^{\Diamond M}_A\!,$
\\[2.5ex]
\>$(\delta^{\Box M})$\hspace{1.7em}\=$\Box\delta^{\Box
M}_A\cirk\delta^{\Box M}_A=\delta^{\Box\Box}_{MA}\cirk\delta^{\Box
M}_A$,\hspace{2em}\= $(\delta^{\Diamond M})$\hspace{1.7em}\=
$\delta^{\Diamond M}_A\cirk\Diamond\delta^{\Diamond
M}_A=\delta^{\Diamond M}_A\cirk \delta^{\Diamond\Diamond}_{MA}$,
\\[1.5ex]
\>$(\Box M\beta)$\> $\varepsilon^\Box_{MA}\cirk\delta^{\Box
M}_A=\mj_{MA}$,\> $(\Diamond M\beta)$\>$\delta^{\Diamond M}_A\cirk
\varepsilon^\Diamond_{MA}=\mj_{MA}$,
\\[1.5ex]
\> ($\delta${\lat N})\> $\Box\delta^{\Diamond
M}_A\cirk\delta^{\Box\Diamond}_{MA}=\delta^{\Box
M}_A\cirk\delta^{\Diamond M}_A$,\> ($\delta${\cir
I})\>$\delta^{\Diamond\Box}_{MA}\cirk \Diamond\delta^{\Box
M}_A=\delta^{\Box M}_A\cirk\delta^{\Diamond M}_A$.
\end{tabbing}

The equations ($\delta^{\Box M}$~{\it nat}), $(\delta^{\Box M})$
and $(\Box M\beta)$ for $M$ being $\Box$ were already assumed for
$S4_\Box$ and $S4_{\Box\Diamond}$, while the equations
($\delta^{\Diamond M}$~{\it nat}), $(\delta^{\Diamond M})$ and
$(\Diamond M\beta)$ for $M$ being $\Diamond$ were already assumed
for $S4_\Diamond$ and $S4_{\Box\Diamond}$. There is no
generalization with $M$ of the equation $(\Box\Box\eta)$ of
Section~3, and of the dual equation for $\Diamond$. These
equations are assumed for $S5_{\Box\Diamond}$ as they were assumed
before for $S4_{\Box\Diamond}$.

The names of the equations ($\delta${\lat N}) and ($\delta${\cir
I}) are derived from graphs related to their left-hand sides (as
will be explained below). These equations are related to the
\emph{Frobenius equations} of Frobenius algebras (see \cite{K03};
for some history concerning the Frobenius equations, see
\cite{K06}, which traces the equations to \cite{CW87}, where they
occur in a different context). The difference is that in the
Frobenius equations $\Box$ and $\Diamond$ are not distinguished.
The equations
\[
\Box\varepsilon^\Box_A\cirk\Box\delta^{\Diamond\Box}_A\cirk\delta^{\Box\Diamond}_{\Box
A}= \varepsilon^\Box_{\Box A}\cirk\delta^{\Diamond\Box}_{\Box
A}\cirk\Diamond\delta^{\Box\Box}_A=\delta^{\Diamond\Box}_A,
\]
or, alternatively, the dual equations
\[
\delta^{\Diamond\Box}_{\Diamond
A}\cirk\Diamond\delta^{\Box\Diamond}_A\cirk\Diamond\varepsilon^\Diamond_A=
\Box\delta^{\Diamond\Diamond}_A\cirk\delta^{\Box\Diamond}_{\Diamond
A}\cirk\varepsilon^\Diamond_{\Diamond A}=\delta^{\Box\Diamond}_A,
\]
suggested by Lawvere (see \cite{LAW69}, p.\ 152, where $\Box$ and
$\Diamond$ are not distinguished), could replace the equations
($\delta${\lat N}) and ($\delta${\cir I}) in our axiomatization of
the equations of $S5_{\Box\Diamond}$.

The equations $(\delta^{\Box M})$ and $(\delta^{\Diamond M})$ are
redundant in this axiomatization. For $(\delta^{\Box M})$ we have
\begin{tabbing}
\hspace{4em}$\Box\delta^{\Box M}_A\cirk\delta^{\Box
M}_A$\=$=\Box\delta^{\Box M}_A\cirk\delta^{\Box
M}_A\cirk\delta^{\Diamond M}_A\cirk\varepsilon^\Diamond_{MA}$,
\hspace{.5em}with $(\Diamond M\beta)$,
\\*[1ex]
\>$=\Box\delta^{\Box M}_A\cirk\Box\delta^{\Diamond
M}_A\cirk\delta^{\Box\Diamond}_{MA}\cirk\varepsilon^\Diamond_{MA}$,
\hspace{.5em}with ($\delta${\lat N}),
\\[1ex]
\>$=\Box\delta^{\Diamond\Box}_{MA}\cirk\Box\Diamond\delta^{\Box
M}_A\cirk\delta^{\Box\Diamond}_{MA}\cirk\varepsilon^\Diamond_{MA}$,
\hspace{.5em}with ($\delta${\cir I}),
\\[1ex]
\>$=\Box\delta^{\Diamond\Box}_{MA}\cirk\delta^{\Box\Diamond}_{\Box
MA}\cirk\varepsilon^\Diamond_{\Box MA}\cirk\delta^{\Box M}_A\!$,
\hspace{.5em}with naturality equations,
\\[1ex]
\>$=\delta^{\Box\Box}_{MA}\cirk\delta^{\Diamond\Box}_{MA}\cirk\varepsilon^\Diamond_{\Box
MA}\cirk\delta^{\Box M}_A\!$, \hspace{.5em}with ($\delta${\lat
N}),
\\*[1ex]
\>$=\delta^{\Box\Box}_{MA}\cirk\delta^{\Box M}_A$,
\hspace{.5em}with $(\Diamond M\beta)$,
\end{tabbing}
and we proceed analogously for $(\delta^{\Diamond M})$ (for an
analogous derivation see \cite{K03}, Proposition 2.3.24, which in
Section 2.3.25 is credited to \cite{Q95}). The equations
$(\delta^{\Box M})$ and $(\delta^{\Diamond M})$ do not seem
however to be redundant if we replace ($\delta${\lat N}) and
($\delta${\cir I}) by the equations suggested by Lawvere.

For $S5_{\Box\Diamond}$, we derive from ($\delta${\lat N}) and
$(\Diamond\Box\beta)$ the equation
\[
\delta^{\Box\Box}_A=\Box\delta^{\Diamond\Box}_A\cirk\delta^{\Box\Diamond}_{\Box
A}\cirk\varepsilon^\Diamond_{\Box A},
\]
and we derive analogously from ($\delta${\cir I}) and
$(\Box\Diamond\beta)$ the equation
\[
\delta^{\Diamond\Diamond}_A=\varepsilon^{\Box}_{\Diamond
A}\cirk\delta^{\Diamond\Box}_{\Diamond
A}\cirk\Diamond\delta^{\Box\Diamond}_A,
\]
which means that the arrows $\delta^{\Box\Box}_A$ and
$\delta^{\Diamond\Diamond}_A$ may be defined in terms of other
arrows, and need not be taken as primitive.

We will now define a category called \emph{Gen}, which will
replace \emph{Rel} to define a functor $G$ from
$S5_{\Box\Diamond}$. The objects of \emph{Gen} are again the
finite ordinals. An arrow of \emph{Gen} from $n$ to $m$ is an
equivalence relation defined on the disjoint union of $n$ and $m$,
which is called a \emph{split equivalence}. The identity arrow
from $n$ to $n$ is the split equivalence that corresponds to the
following picture:
\begin{center}
\begin{picture}(32,40)

\put(0,10){\circle{2}} \put(32,10){\circle{2}}
\put(0,30){\circle{2}} \put(32,30){\circle{2}}

\put(10,20){$\ldots$}

\put(0,12){\line(0,1){16}} \put(32,12){\line(0,1){16}}

\put(0,7){\makebox(0,0)[t]{\scriptsize$n\mn 1$}}
\put(32,7){\makebox(0,0)[t]{\scriptsize$0$}}
\put(0,33){\makebox(0,0)[b]{\scriptsize$n\mn 1$}}
\put(32,33){\makebox(0,0)[b]{\scriptsize$0$}}

\end{picture}
\end{center}
which is empty if $n=0$. We do not draw in such pictures the loops
corresponding to the pairs ${(x,x)}$. Composition of arrows is
defined, roughly speaking, as the transitive closure of the union
of the two relations composed, where we omit the ordered pairs one
of whose members is in the middle (see \cite{DP03a}, Section~2,
and \cite{DP03b}, Section~2, for a detailed definition). For
example, the split equivalences $R_1$ and $R_2$ corresponding to
the following two pictures
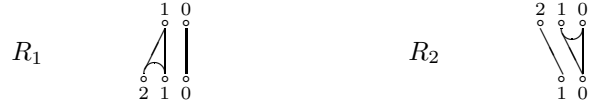
\begin{center}
\begin{picture}(264,40)

\put(50,10){\circle{2}} \put(58,10){\circle{2}}
\put(66,10){\circle{2}} \put(58,31){\circle{2}}
\put(66,31){\circle{2}}

\put(58,29){\line(-1,-2){7.6}} \put(58,12){\line(0,1){17}}
\put(66,12){\line(0,1){17}}

\put(54,12){\oval(8,8)[t]}

\put(50,7){\makebox(0,0)[t]{\scriptsize$2$}}
\put(58,7){\makebox(0,0)[t]{\scriptsize$1$}}
\put(66,7){\makebox(0,0)[t]{\scriptsize$0$}}
\put(58,34){\makebox(0,0)[b]{\scriptsize$1$}}
\put(66,34){\makebox(0,0)[b]{\scriptsize$0$}}

\put(0,20){\makebox(0,0)[l]{$R_1$}}

\put(200,31){\circle{2}} \put(208,10){\circle{2}}
\put(216,10){\circle{2}} \put(208,31){\circle{2}}
\put(216,31){\circle{2}}

\put(216,12){\line(-1,2){7.6}} \put(216,12){\line(0,1){17}}
\put(208.5,12){\line(-1,2){8.5}}

\put(212,29){\oval(8,8)[b]}

\put(200,34){\makebox(0,0)[b]{\scriptsize$2$}}
\put(208,7){\makebox(0,0)[t]{\scriptsize$1$}}
\put(216,7){\makebox(0,0)[t]{\scriptsize$0$}}
\put(208,34){\makebox(0,0)[b]{\scriptsize$1$}}
\put(216,34){\makebox(0,0)[b]{\scriptsize$0$}}

\put(150,20){\makebox(0,0)[l]{$R_2$}}

\end{picture}
\end{center}
are composed as follows, so as to yield the split equivalence
${R_2\cirk R_1}$ that corresponds to the picture on the right-hand
side
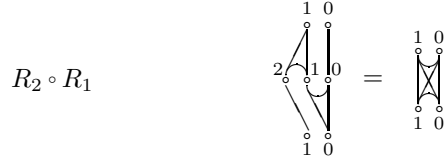
\begin{center}
\begin{picture}(164,60)

\put(50,30){\circle{2}} \put(58,30){\circle{2}}
\put(66,30){\circle{2}} \put(58,9){\circle{2}}
\put(66,9){\circle{2}} \put(58,51){\circle{2}}
\put(66,51){\circle{2}}

\put(58,49){\line(-1,-2){7.6}} \put(58,32){\line(0,1){17}}
\put(66,32){\line(0,1){17}}

\put(54,32){\oval(8,8)[t]}

\put(66,11){\line(-1,2){7.6}} \put(66,11){\line(0,1){17}}
\put(58.5,11){\line(-1,2){8.5}}

\put(62,28){\oval(8,8)[b]}

\put(58,6){\makebox(0,0)[t]{\scriptsize$1$}}
\put(66,6){\makebox(0,0)[t]{\scriptsize$0$}}
\put(58,55){\makebox(0,0)[b]{\scriptsize$1$}}
\put(66,55){\makebox(0,0)[b]{\scriptsize$0$}}
\put(49,32){\makebox(0,0)[br]{\scriptsize$2$}}
\put(59,32){\makebox(0,0)[bl]{\scriptsize$1$}}
\put(67,32){\makebox(0,0)[bl]{\scriptsize$0$}}

\put(-54,30){\makebox(0,0)[l]{$R_2\cirk R_1$}}

\put(100,19){\circle{2}} \put(108,19){\circle{2}}
\put(100,41){\circle{2}} \put(108,41){\circle{2}}
\put(104,21){\oval(8,8)[t]} \put(104,39){\oval(8,8)[b]}
\put(100.8,23.4){\line(1,2){7}} \put(107.2,23.4){\line(-1,2){7}}
\put(100,21){\line(0,1){18}} \put(108,21){\line(0,1){18}}
\put(100,16){\makebox(0,0)[t]{\scriptsize$1$}}
\put(108,16){\makebox(0,0)[t]{\scriptsize$0$}}
\put(100,44){\makebox(0,0)[b]{\scriptsize$1$}}
\put(108,44){\makebox(0,0)[b]{\scriptsize$0$}}

\put(83,30){\makebox(0,0){$=$}}

\end{picture}
\end{center}

We define the functor $G$ from $S5_{\Box\Diamond}$ to \emph{Gen}
by stipulating first that $GA$ is the length of the object $A$. On
arrows, we have first that $G\mj_A$ is the identity arrow of
\emph{Gen} from $GA$ to $GA$. For ${GA=n}$, let
$G\varepsilon^\Box_A$ and $G\varepsilon^\Diamond_A$ be the split
equivalences that correspond respectively to the pictures given
for $G\varepsilon^\Box_n$ and $G\varepsilon^\Diamond_n$ in
Section~3. We have next, for ${GA=n}$, the clauses that correspond
to the following pictures:
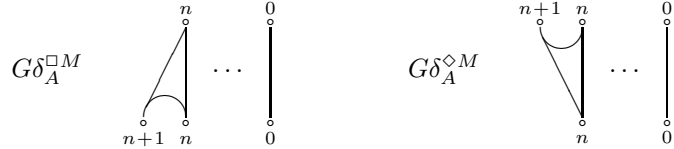
\begin{center}
\begin{picture}(277,56)

\put(50,10){\circle{2}} \put(66,10){\circle{2}}
\put(98,10){\circle{2}} \put(66,48){\circle{2}}
\put(98,48){\circle{2}}

\put(76,29){$\ldots$}

\put(66,46){\line(-1,-2){15}} \put(66,12){\line(0,1){34}}
\put(98,12){\line(0,1){34}}

\put(58,12){\oval(16,16)[t]}

\put(50,2){\makebox(0,0)[b]{\scriptsize$n\pl 1$}}
\put(66,2){\makebox(0,0)[b]{\scriptsize$n$}}
\put(98,2){\makebox(0,0)[b]{\scriptsize$0$}}
\put(66,51){\makebox(0,0)[b]{\scriptsize$n$}}
\put(98,51){\makebox(0,0)[b]{\scriptsize$0$}}

\put(0,30){\makebox(0,0)[l]{$G\delta^{\Box M}_A$}}

\put(200,48){\circle{2}} \put(216,10){\circle{2}}
\put(248,10){\circle{2}} \put(216,48){\circle{2}}
\put(248,48){\circle{2}}

\put(226,29){$\ldots$}

\put(216,12){\line(-1,2){15}} \put(216,12){\line(0,1){34}}
\put(248,12){\line(0,1){34}}

\put(208,46){\oval(16,16)[b]}

\put(200,51){\makebox(0,0)[b]{\scriptsize$n\pl 1$}}
\put(216,2){\makebox(0,0)[b]{\scriptsize$n$}}
\put(248,2){\makebox(0,0)[b]{\scriptsize$0$}}
\put(216,51){\makebox(0,0)[b]{\scriptsize$n$}}
\put(248,51){\makebox(0,0)[b]{\scriptsize$0$}}

\put(150,30){\makebox(0,0)[l]{$G\delta^{\Diamond M}_A$}}

\end{picture}
\end{center}
The semicircle joining $n$ and ${n\pl 1}$ at the bottom (in the
target) in the left picture is the \emph{cap} ${(n,n\pl 1)}$, and
the semicircle joining $n$ and ${n\pl 1}$ at the top (in the
source) in the right picture is the \emph{cup} ${(n,n\pl 1)}$.
These two pictures are like those we had in Section~3 for
$G\delta^{\Box\Box}_n$ and $G\delta^{\Diamond\Diamond}_n$ but with
the cap and the cup added.

As before, we have ${G(g\cirk f)=Gg\cirk Gf}$, and for
${Gf\!:n\vdash m}$ the partition induced by the split equivalence
$GMf$ is obtained from the partition induced by the split
equivalence $Gf$ by adding the equivalence class ${\{n,m\}}$,
where $n$ is in the source and $m$ in the target. We easily check
by induction on the length of derivation that if ${f=g}$ in
$S5_{\Box\Diamond}$, then ${Gf=Gg}$ in \emph{Gen}; hence $G$ so
defined is indeed a functor.

The split equivalences $R_1$ and $R_2$ in the example above may be
taken to be respectively $G\delta^{\Box\Diamond}_{MA}$ and
$G\Box\delta^{\Diamond M}_A$ for $A$ being empty. Then ${R_2\cirk
R_1}$ is the $G$ image of an instance of the left-hand side of
($\delta${\lat N}), and when in the left picture corresponding to
${R_2\cirk R_1}$ we omit the cup ${(0,1)}$ and the cap ${(1,2)}$
in the middle, we obtain the form of {\lat N}. (This explains
{\lat N} in the name of ($\delta${\lat N}); horizontally, we would
obtain {\lat Z}, and in the comments in \cite{K06} this horizontal
look at the matter is favoured. The {\cir I} of the name of
($\delta${\cir I}) arises analogously.)

Before proving that this functor is faithful, note that the
coherence results established in the preceding text with respect
to \emph{Rel} could be established with respect to \emph{Gen}, by
relying on functors $G$ obtained by restricting appropriately the
functor $G$ from $S5_{\Box\Diamond}$ to \emph{Gen}. For that we
have to check first that these restricted functors are indeed
functors, which is done by induction on the length of derivation
(the essential ingredient in this induction is to go through the
axiomatic equations). This is nearly all we have to check, because
the faithfulness of these functors can next be established by
proceeding as before, via the same normal forms. Roughly speaking,
adding the cups and caps to the pictures we had before does not
change matters. (For a more detailed treatment of the relationship
between \emph{Rel} and \emph{Gen} see \cite{DP09}.)

Next, as an auxiliary result, we establish coherence with respect
to \emph{Gen} for the category $S5_{\Box\Diamond}^\str$, defined
by omitting from the definition of $S5_{\Box\Diamond}$ the arrow
terms $\varepsilon_A^\Diamond$ and $\delta_A^{\Box M}$, and all
the equations involving them explicitly. This means that we have
in $S5_{\Box\Diamond}^\str$ only the primitive arrow terms
$\mj_A$, $\varepsilon_A^\Box$ and $\delta_A^{\Diamond M}$, for
which we assume the categorial and functorial equations plus
($\varepsilon^\Box$~{\it nat}), ($\delta^{\Diamond M}$~{\it nat})
and $(\delta^{\Diamond M})$. The functor $G$ from
$S5_{\Box\Diamond}^\str$ to \emph{Gen} is defined by restricting
the definition of $G$ from $S5_{\Box\Diamond}$ to \emph{Gen}. We
have the following.

\prop{$S5_{\Box\Diamond}^\str$ Coherence}{The functor $G$ from
$S5_{\Box\Diamond}^\str$ to {\it Gen} is faithful.}

\dkz Suppose that for $f$ and $g$ arrow terms of
$S5_{\Box\Diamond}^\str$ of the same type we have ${Gf=Gg}$. We
prove that ${f=g}$ in $S5_{\Box\Diamond}^\str$ by induction on the
number $n$ of occurrences of $\delta^{\Diamond M}$ in $f$, which
must be equal to that number for $g$. If ${n=0}$, then we rely on
$S$ Coherence of Section~2. If ${n>0}$, then we rely on a lemma
that says that if in the picture corresponding to $Gf$ we have a
cup ${(i,i\pl 1)}$ in the source, then $f$ is equal in
$S5_{\Box\Diamond}^\str$ to an arrow term of the form ${f'\cirk
A\delta^{\Diamond M}_B}$ such that ${GB=i}$. This lemma is
sufficient because if there are no cups ${(i,i\pl 1)}$ in the
source, then $f$ and $g$ are equal respectively to ${f'\cirk h}$
and ${g'\cirk h}$ for $h$ without $\delta^{\Diamond M}$, and a cup
${(i,i\pl 1)}$ in the source of $Gf'$, which is equal to $Gg'$.

Here is a sketch of the proof of this lemma. We first transform
$f$ into the developed form ${f_n\cirk\ldots\cirk f_1}$ (see
Section~2), and then we find the $f_i$ ``responsible'' for the cup
${(i,i\pl 1)}$. We use then the equations of
$S5_{\Box\Diamond}^\str$, and we may rely in particular on
$(\delta^{\Diamond M})$, to permute this $f_i$ to the right, until
a descendent of it becomes the rightmost factor. \qed

\vspace{2ex}

Let $S5_{\Box\Diamond}^\rts$ be the category isomorphic to
$(S5_{\Box\Diamond}^\str)^{op}$, where $\varepsilon^\Box$ and
$\delta^{\Diamond M}$ are replaced by $\varepsilon^\Diamond$ and
$\delta^{\Box M}$ respectively. The equations for
$S5_{\Box\Diamond}^\rts$ are dual to those for
$S5_{\Box\Diamond}^\str$ (instead of ($\varepsilon^\Box$~{\it
nat}), ($\delta^{\Diamond M}$~{\it nat}) and $(\delta^{\Diamond
M})$ we have ($\varepsilon^\Diamond$~{\it nat}), ($\delta^{\Box
M}$~{\it nat}) and $(\delta^{\Box M})$). Coherence for
$S5_{\Box\Diamond}^\str$, which we have established above,
delivers of course coherence for $S5_{\Box\Diamond}^\rts$. We can
then establish the following.

\prop{$S5_{\Box\Diamond}$ Coherence}{The functor $G$ from
$S5_{\Box\Diamond}$ to {\it Gen} is faithful.}

\dkz We verify that every arrow term $f$ of $S5_{\Box\Diamond}$ is
equal in $S5_{\Box\Diamond}$ to an arrow term in the \emph{normal
form} ${f_2\cirk f_1}$ where $f_1$ is an arrow term of
$S5_{\Box\Diamond}^\str$ and $f_2$ is an arrow term of
$S5_{\Box\Diamond}^\rts$. It is easy to see that $Gf$ determines
uniquely $Gf_1$ and $Gf_2$, as well as the target of $f_1$. To
conclude the proof of $S5_{\Box\Diamond}$ Coherence we rely then
on coherence for $S5_{\Box\Diamond}^\str$ and
$S5_{\Box\Diamond}^\rts$.\qed

\vspace{2ex}

\noindent The normal form we have used in this proof is of the
thin kind (cf.\ Section~4).

Suppose that in the definition of $G$ for $S5_{\Box\Diamond}$ we
omit from the picture corresponding to the clause for
$G\delta_A^{\Box M}$ the cap ${(n,n\pl 1)}$, and from the picture
corresponding to the clause for $G\delta_A^{\Diamond M}$ the cup
${(n,n\pl 1)}$. The target category for that $G$ would be
\emph{Rel}, but we could not show that this defines a functor from
$S5_{\Box\Diamond}$, because of the equations ($\delta${\lat N})
and ($\delta${\cir I}). These equations require the caps and cups,
and the split equivalences of \emph{Gen}.

We can prove coherence for $S5_{\Box\Diamond}$ with respect to a
functor $G^d$ from $S5_{\Box\Diamond}$ to \emph{Gen} that is a
kind of dual of the functor $G$ we had above. It interchanges the
role of $\varepsilon$ and $\delta$ in the following manner. On
objects, $G^dA$ is ${GA\pl 1}$. On arrows, we have
\begin{tabbing}
\hspace{7em}\=$\;\;\;G^d\varepsilon^\Box_A$\=$=G\delta^{\Diamond\Diamond}_A$,
\hspace{5em}\=$\;\;\;G^d\varepsilon^\Diamond_A$\=$=G\delta^{\Box\Box}_A$,
\\*[1.5ex]
\>$G^d\delta^{\Box
M}_A$\>$=GM\varepsilon^\Diamond_{MA}$,\>$G^d\delta^{\Diamond
M}_A$\>$=GM\varepsilon^\Box_{MA}$,
\\[2ex]
\hspace{12em}$G^d(g\cirk f)=G^dg\cirk G^df$;
\end{tabbing}
for ${G^dMf}$ we have a clause exactly analogous to the clause for
${GMf}$ for $S5_{\Box\Diamond}$. Graphically, for the length of
$A$ being $n$, we have the following:
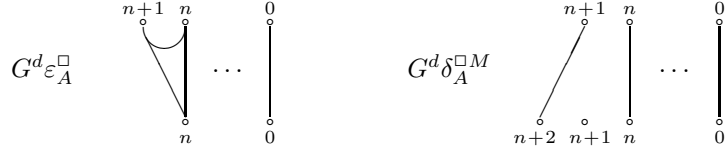
\begin{center}
\begin{picture}(278,56)

\put(200,10){\circle{2}} \put(217,10){\circle{2}}
\put(234,10){\circle{2}} \put(268,10){\circle{2}}
\put(217,48){\circle{2}} \put(234,48){\circle{2}}
\put(268,48){\circle{2}}

\put(245,29){$\ldots$}

\put(200,12){\line(1,2){17}} \put(234,12){\line(0,1){34}}
\put(268,12){\line(0,1){34}}

\put(198,2){\makebox(0,0)[b]{\scriptsize$n\pl 2$}}
\put(219,2){\makebox(0,0)[b]{\scriptsize$n\pl 1$}}
\put(234,2){\makebox(0,0)[b]{\scriptsize$n$}}
\put(268,2){\makebox(0,0)[b]{\scriptsize$0$}}
\put(217,51){\makebox(0,0)[b]{\scriptsize$n\pl 1$}}
\put(234,51){\makebox(0,0)[b]{\scriptsize$n$}}
\put(268,51){\makebox(0,0)[b]{\scriptsize$0$}}

\put(150,30){\makebox(0,0)[l]{$G^d\delta^{\Box M}_A$}}

\put(50,48){\circle{2}} \put(66,10){\circle{2}}
\put(98,10){\circle{2}} \put(66,48){\circle{2}}
\put(98,48){\circle{2}}

\put(76,29){$\ldots$}

\put(66,12){\line(-1,2){15}} \put(66,12){\line(0,1){34}}
\put(98,12){\line(0,1){34}}

\put(58,46){\oval(16,16)[b]}

\put(50,51){\makebox(0,0)[b]{\scriptsize$n\pl 1$}}
\put(66,2){\makebox(0,0)[b]{\scriptsize$n$}}
\put(98,2){\makebox(0,0)[b]{\scriptsize$0$}}
\put(66,51){\makebox(0,0)[b]{\scriptsize$n$}}
\put(98,51){\makebox(0,0)[b]{\scriptsize$0$}}

\put(0,30){\makebox(0,0)[l]{$G^d\varepsilon^{\Box}_A$}}

\end{picture}
\end{center}
and analogously for $G^d\varepsilon^\Diamond_A$ and
$G^d\delta^{\Diamond M}_A$.

That $G^d$ is indeed a functor is checked by induction on the
length of derivation of the equations of $S5_{\Box\Diamond}$. The
only problematic case arises with the equations
($\varepsilon^M$~{\it nat}), where we rely on the fact that the
pair ${(n\mn 1,m\mn 1)}$ belongs to ${G^df\!:n\vdash m}$. That
$G^d$ is a faithful functor can be shown either directly, as
$S5_{\Box\Diamond}$ Coherence above, via the same normal form, or,
alternatively, we can rely on the Maximality of
$S5_{\Box\Diamond}$ of Section 11 (which presupposes
$S5_{\Box\Diamond}$ Coherence).

This duality between $\varepsilon$ and $\delta$, exhibited by
$G^d$, was already present in the category $S_+$ of Section~2,
whose arrows $\xi_{n+1}$ could be interpreted either as
$\varepsilon_{n+1}$ or as $\delta_n$ arrows. Functors dual to the
functors $G$ from $S4_\Box$, $S4_\Diamond$ and $S4_{\Box\Diamond}$
to \emph{Rel}, as $G^d$ is dual to $G$ from $S5_{\Box\Diamond}$ to
\emph{Gen}, can be defined analogously (just omit the cups and
caps from the $G^d$ images). The faithfulness of these dual
functors can be proved either directly, via normal forms used
previously, or for $S4_\Box$ and $S4_\Diamond$ we could rely on
their maximality (see Section~9). We could also rely on a result
about the duality of the simplicial category, analogous to the
duality between $G$ and $G^d$, which is explained in \cite{DP08a}
(end of Section~6).

Note that in $S5_{\Box\Diamond}$ we have the arrows
\[
\Box\varepsilon^\Diamond_A\cirk\delta^{\Diamond\Box}_A\!:\Diamond\Box
A\vdash\Box\Diamond A,
\hspace{5em}\delta^{\Box\Diamond}_A\cirk\Diamond\varepsilon^\Box_A\!:\Diamond\Box
A\vdash\Box\Diamond A,
\]
which are of the same type as the arrows $\chi^{\Diamond\Box}_A$
of $S4.2_{\Box\Diamond}$ (see the preceding section), but the
equation $(\varepsilon^\Box\chi^{\Diamond\Box})$ fails for the
first arrow, and the equation
$(\varepsilon^\Diamond\chi^{\Diamond\Box})$ fails for the second,
as can be easily verified with the help of the functor $G$ form
$S5_{\Box\Diamond}$ to \emph{Gen}. That these arrows of
$S5_{\Box\Diamond}$ do not amount to $\chi^{\Diamond\Box}_A$ is
clear from their interpretation via $G$. So, although, as far as
theorems and provable sequents are concerned, the modal logic
$S4.2$ is included in the modal logic $S5$, from a
proof-theoretical point of view we should not assume that $S4.2$
is a subsystem of $S5$. Our $S5_{\Box\Diamond}$ does not cover
$S4.2_{\Box\Diamond}$. There are deductions in $S4.2$ (i.e.\
arrows of $S4.2_{\Box\Diamond}$) absent from $S5$.

In $S5_{\Box\Diamond}$ the endofunctor $\Diamond$ is left adjoint
to the endofunctor $\Box$ (for the notion of adjunction, see the
beginning of Section 10). The members of the unit and counit of
this adjunction are respectively the arrows
\[
\delta^{\Box\Diamond}_A\cirk\varepsilon^\Diamond_A\!:A\vdash\Box\Diamond
A,\hspace{5em}
\varepsilon^\Box_A\cirk\delta^{\Diamond\Box}_A\!:\Diamond\Box
A\vdash A,
\]
which correspond to modal laws found in the modal system $B$ (see
\cite{HC96}, p.\ 62). We will treat of matters pertaining to this
adjunction in Section 10.

\section{The category $5S_{\Box\Diamond}$}

We consider now a category isomorphic to $S5_{\Box\Diamond}$, a
kind of mirror image of it. We define this category like
$S4_{\Box\Diamond}$, save that instead of $\delta$ we write
$\sigma$, and we have the following additions. We have the
additional primitive arrow terms
\[
\sigma^{\Diamond\Box}_A\!:\Diamond A\vdash\Diamond\Box
A,\hspace{5em}\sigma^{\Box\Diamond}_A\!:\Box\Diamond A\vdash\Box
A.
\]
The modal laws corresponding to the types of these arrow terms
were investigated in \cite{LS77} (p.\ 67).

The equations of $5S_{\Box\Diamond}$ are obtained by assuming
those assumed for $S4_{\Box\Diamond}$, with $\delta$ replaced by
$\sigma$, and the following additional equations:
\begin{tabbing}
\hspace{.4em}\=$(\sigma^{M\Box})$\hspace{1.4em}\=$\sigma^{M\Box}_{\Box
A}\cirk\sigma^{M\Box}_A=M\sigma^{\Box\Box}_A\cirk\sigma^{M\Box}_A$,\hspace{1.2em}\=
$(\sigma^{M\Diamond})$\hspace{1.4em}\=
$\sigma^{M\Diamond}_A\cirk\sigma^{M\Diamond}_{\Diamond
A}=\sigma^{M\Diamond}_A\cirk M\sigma^{\Diamond\Diamond}_A$,\kill

\hspace{.4em}($\sigma^{M\Box}$~{\it nat})\hspace{.5em}$M\Box
f\cirk\sigma^{M\Box}_A=\sigma^{M\Box}_B\cirk Mf$,\hspace{.7em}
($\sigma^{M\Diamond}$~{\it nat})\hspace{.5em}
$\sigma^{M\Diamond}_B\cirk M\Diamond
f=Mf\cirk\sigma^{M\Diamond}_A\!,$
\\[2.5ex]
\>$(\Box M\eta)$\>
$M\varepsilon^\Box_A\cirk\sigma^{M\Box}_A=\mj_{MA}$,\> $(\Diamond
M\eta)$\>$\sigma^{M\Diamond}_A\cirk
M\varepsilon^\Diamond_A=\mj_{MA}$,
\\[1.5ex]
\> ($\sigma${\lat N})\>
$M\sigma^{\Box\Diamond}_A\cirk\sigma^{M\Box}_{\Diamond
A}=\sigma^{M\Box}_A\cirk\sigma^{M\Diamond}_A$,\> ($\sigma${\cir
I})\>$\sigma^{M\Diamond}_{\Box A}\cirk
M\sigma^{\Diamond\Box}_A=\sigma^{M\Box}_A\cirk\sigma^{M\Diamond}_A$.
\end{tabbing}
The following equations can be derived (see the derivation of
$(\delta^{\Box M})$ in the preceding section):
\begin{tabbing}
\hspace{.4em}\=$(\sigma^{M\Box})$\hspace{1.4em}\=$\sigma^{M\Box}_{\Box
A}\cirk\sigma^{M\Box}_A=M\sigma^{\Box\Box}_A\cirk\sigma^{M\Box}_A$,\hspace{1.2em}\=
$(\sigma^{M\Diamond})$\hspace{1.4em}\=
$\sigma^{M\Diamond}_A\cirk\sigma^{M\Diamond}_{\Diamond
A}=\sigma^{M\Diamond}_A\cirk M\sigma^{\Diamond\Diamond}_A$.
\end{tabbing}

It is not difficult to show that the categories
$S5_{\Box\Diamond}$ and $5S_{\Box\Diamond}$ are isomorphic. In
this isomorphism, the object $A$ is mapped to $A$ read from right
to left. (This isomorphism does not preserve the functors $\Box$
and $\Diamond$.)

It follows that for $5S_{\Box\Diamond}$ we can establish coherence
with respect to the functor $G$ from $5S_{\Box\Diamond}$ to
\emph{Gen} defined like $G$ from $S5_{\Box\Diamond}$ to
\emph{Gen}; namely, $G\sigma^{M_1M_2}_A=G\delta^{M_2M_1}_A$. In
$5S_{\Box\Diamond}$, the endofunctor $\Diamond$ is right adjoint
to the endofunctor $\Box$, while in $S5_{\Box\Diamond}$ it was
left adjoint, as we noted at the end of the preceding section.

Note that in $5S_{\Box\Diamond}$ we do not have an arrow of the
type ${\emptyset\vdash \Box\emptyset}$ for $\emptyset$ being the
empty sequence. Analogously, we do not have an arrow of the type
${\Diamond\emptyset\vdash\emptyset}$. This is because for every
arrow $f$ of $5S_{\Box\Diamond}$, every occurrence of $\Box$ in
the target of $f$ must be linked by $Gf$ to an occurrence of
$\Box$ in the source of $f$ or an occurrence of $\Diamond$ in the
target of $f$, and every occurrence of $\Diamond$ in the source of
$f$ must be linked by $Gf$ to an occurrence of $\Diamond$ in the
target of $f$ or an occurrence of $\Box$ in the source of $f$.
Another way to conclude that arrows of the type ${\emptyset\vdash
\Box\emptyset}$ or ${\Diamond\emptyset\vdash\emptyset}$ do not
exist in $5S_{\Box\Diamond}$ is to appeal to the isomorphism of
$5S_{\Box\Diamond}$ with $S5_{\Box\Diamond}$, and the well-known
fact that in the modal logic $S5$ we do not have modal laws
corresponding to these types.  However, in the extension of the
modal logic $T$ (namely, the normal modal logic with the axiom
${\Box p\str p}$ or ${p\str\Diamond p}$) with the axiom ${\Diamond
p\str\Diamond\Box p}$ or ${\Box\Diamond p\str\Box p}$, we can
derive ${p\str\Box p}$ and ${\Diamond p\str p}$. We have
\begin{tabbing}
\hspace{9em}\= $\Box p\str\Diamond\Box p$, \hspace{1em}\= by
$\alpha\str\Diamond\alpha$,
\\*[.5ex]
\> $\Diamond(p\str\Box p)$,\> by laws of normal modal logics,
\\[.5ex]
\> $\Box\Diamond(p\str\Box p)$,\> by necessitation,
\\[.5ex]
\> $\Box(p\str\Box p)$,\> by $\Box\Diamond\alpha\str\Box\alpha$,
\\*[.5ex]
\> $p\str\Box p$,\> by $\Box\alpha\str\alpha$.
\end{tabbing}
This may be the reason why the modalities of $5S_{\Box\Diamond}$
are not usually considered, though the laws governing these
modalities are as interesting as those of $S5_{\Box\Diamond}$,
whose faithful image they are.

\section{Trijunctions, dyads and codyads}

In this section we show that the assumptions made for the category
$S5_{\Box\Diamond}$ can be justified by adjunctions underlying the
comonad and monad structures of that category.

A \emph{comonad} on a category $\cal C$ is a structure
${\langle{\cal C},\Box,\varepsilon^\Box,\delta^{\Box\Box}\rangle}$
where $\Box$ is an endofunctor of $\cal C$, while
${\varepsilon^\Box\!:\Box\strt I_{\cal C}}$ and
${\delta^{\Box\Box}\!:\Box\strt\Box\Box}$, for $I_{\cal C}$ being
the identity functor of $\cal C$, are natural transformations that
satisfy the equations of $S4_\Box$ (provided $n$ and ${n\pl 1}$
are replaced respectively by $A$ and $\Box A$, for $A$ an object
of $\cal C$). The category $S4_\Box$ is the free comonad generated
by a single object (understood as an arrowless one-node graph, or
the trivial one-object category; for details, see \cite{D99},
Chapter~5, and \cite{D08}, Section~4). A \emph{monad} on $\cal C$
is a structure ${\langle{\cal
C},\Diamond,\varepsilon^\Diamond,\delta^{\Diamond\Diamond}\rangle}$
defined analogously by reference to $S4_\Diamond$, which is the
free monad generated by a single object.

We call \emph{dyad} on $\cal C$ a structure that includes a
comonad on $\cal C$, a monad on $\cal C$, and two additional
natural transformations
${\delta^{\Box\Diamond}\!:\Diamond\strt\Box\Diamond}$ and
${\delta^{\Diamond\Box}\!:\Diamond\Box\strt\Box}$ that satisfy the
equations of $S5_{\Box\Diamond}$. The category $S5_{\Box\Diamond}$
is the free dyad generated by a single object.

We call \emph{codyad} on $\cal C$ a structure that includes a
comonad on $\cal C$, a monad on $\cal C$, and two additional
natural transformations
${\sigma^{\Diamond\Box}\!:\Diamond\strt\Diamond\Box}$ and
${\sigma^{\Box\Diamond}\!:\Box\Diamond\strt\Box}$ that satisfy the
equations of $5S_{\Box\Diamond}$. The category $5S_{\Box\Diamond}$
is the free codyad generated by a single object.

A \emph{trijunction} is a structure made of the categories $\cal
A$ and $\cal B$, the functor $U$ from $\cal A$ to $\cal B$, and
the functors $L$ and $R$ from $\cal B$ to $\cal A$, such that $L$
is left adjoint to $U$, with the counit $\varphi^L\!:LU\strt
I_{\cal A}$ and unit $\gamma^L\!:I_{\cal B}\strt UL$, and $R$ is
right adjoint to $U$, with the counit $\varphi^R\!:UR\strt I_{\cal
B}$ and unit $\gamma^R\!:I_{\cal A}\strt RU$ (for the notion of
adjunction, see the beginning of Section 10).

The notion of trijunction is very well known, but no special name
seems to be commonly used for it. An important example of a
trijunction is obtained when $\cal A$ is a category with products
and coproducts; then $\cal B$ is the product category ${\cal
A}\times{\cal A}$, the functor $U$ is the diagonal functor, and
the functors $L$ and $R$ are respectively the coproduct and
product bifunctors. Another example of a trijunction, interesting
for logic, which involves the functor of substitution and the
existential and universal quantifiers, may be found in Lawvere's
hyperdoctrines (see \cite{LAW69a} and \cite{LAW70}). A trijunction
involving the category of adjunctions, the category of monads (or
comonads), and the Eilenberg-Moore and the Kleisli constructions
is investigated in \cite{P70} (see also \cite{D99}, Sections
5.2.3-4; cf.\ also \cite{A74}). Trijunctions, and in connection
with them the adjunction from the end of Section~6, are mentioned
in \cite{AW06} (Section 10.4). Particular trijunctions are called
\emph{quasi-Frobenius triples of functors} in \cite{CIN06}. In
\cite{DP08a}, the trijunctions where the functors $L$ and $R$ are
the same functor are called \emph{bijunctions}, and trijunctions
where $U$, $L$ and $R$ are all the same endofunctor are
\emph{self-adjunctions} (examples of such structures may be found
in \cite{DP03c}; see also \cite{DP08a}).

The relationship between the notions of trijunction, dyad and
codyad is analogous to a certain extent to the relationship
between the notions of adjunction, monad and comonad. Every
trijunction gives rise to a dyad on $\cal B$ with $\Box$ being
$UR$ and $\Diamond$ being $UL$; for $B$ an object of $\cal B$, we
have
\begin{tabbing}
\hspace{6.5em}\=$\varepsilon^\Box_B=_{df}\varphi^R_B$,
\hspace{2.25em}\=$\delta^{\Box\Box}_B=_{df}U\gamma^R_{RB}$,
\hspace{2.25em}\=$\delta^{\Box\Diamond}_B=_{df}U\gamma^R_{LB}$,
\\*[1ex]
\>$\varepsilon^\Diamond_B=_{df}\gamma^L_B$,
\>$\delta^{\Diamond\Diamond}_B=_{df}U\varphi^L_{LB}$,
\>$\delta^{\Diamond\Box}_B=_{df}U\varphi^L_{RB}$.
\end{tabbing}
Every trijunction gives analogously rise to a codyad on $\cal A$
with $\Box$ being $LU$ and $\Diamond$ being $RU$. Conversely,
every dyad or codyad gives rise to a trijunction by a construction
analogous to the Eilenberg-Moore construction of an adjunction out
of a monad or comonad (see \cite{ML98}, Sections VI.2, and
\cite{D99}, Sections 5.1.7). We present here this construction.

For a dyad on $\cal C$, let ${\cal C}^\Box_\Diamond$ be the
category whose objects are of the form ${\langle A,d,g\rangle}$
for ${d\!:A\str\Box A}$ and ${g\!:\Diamond A\str A}$ arrows of
$\cal C$ that satisfy the conditions below. Strictly speaking, the
mentioning of the object $A$ is here superfluous, but it is kept
to be in tune with common usage concerning the Eilenberg-Moore
construction. The conditions for $d$ and $g$ are the following
equations, analogous to the similarly named equations of
$S5_{\Box\Diamond}$ in Section~6:
\begin{tabbing}
\hspace{1.7em}\= $(\varepsilon^\Box d)$\hspace{2.5em}\=
$\varepsilon^\Box_A\cirk d=\mj_A$,\hspace{7em}\=
$(\varepsilon^\Diamond g)$\hspace{2.5em}\=
$g\cirk\varepsilon^\Diamond_A=\mj_A$, \hspace{1.7em}\=
$(\varepsilon^\Box d)$\hspace{2.5em}\= $\varepsilon^\Box_A\cirk
d=\mj_A$,\hspace{7em}\= $(\varepsilon^\Diamond g)$\hspace{2.5em}\=
$g\cirk\varepsilon^\Diamond_A=\mj_A$,\kill

\>$(\Box M\beta\;d)$\>$\varepsilon^\Box_A\cirk d=\mj_A$,\>
$(\Diamond M\beta\;g)$\>$g\cirk\varepsilon^\Diamond_A=\mj_A$,
\\[1.5ex]
\> ($\delta${\lat N}$\;g$)\> $\Box
g\cirk\delta^{\Box\Diamond}_A=d\cirk g$,\> ($\delta${\cir
I}$\;d$)\> $\delta^{\Diamond\Box}_A\cirk\Diamond d=d\cirk g$.
\end{tabbing}
The equations
\begin{tabbing}
\hspace{1.7em}\= $(\varepsilon^\Box d)$\hspace{2.5em}\=
$\varepsilon^\Box_A\cirk d=\mj_A$,\hspace{7em}\=
$(\varepsilon^\Diamond g)$\hspace{2.5em}\=
$g\cirk\varepsilon^\Diamond_A=\mj_A$, \hspace{1.7em}\=
$(\varepsilon^\Box d)$\hspace{2.5em}\= $\varepsilon^\Box_A\cirk
d=\mj_A$,\hspace{7em}\= $(\varepsilon^\Diamond g)$\hspace{2.5em}\=
$g\cirk\varepsilon^\Diamond_A=\mj_A$,\kill

\>$(\delta^{\Box M}\,d)$\>$\Box d\cirk d=\delta^{\Box\Box}_A\cirk
d$,\>$(\delta^{\Diamond M}\,g)$\> $g\cirk\Diamond
g=g\cirk\delta^{\Diamond\Diamond}_A$
\end{tabbing}
can be derived (see the derivation of $(\delta^{\Box M})$ in
Section~6). An arrow of ${\cal C}^\Box_\Diamond$ from ${\langle
A_1, d_1,g_1\rangle}$ to ${\langle A_2,d_2,g_2\rangle}$ is an
arrow ${h\!:A_1\str A_2}$ of $\cal C$, indexed by ${\langle A_1,
d_1,g_1\rangle}$ and ${\langle A_2,d_2,g_2\rangle}$, such that the
following equations hold:
\begin{tabbing}
\hspace{1.7em}\= $(\varepsilon^\Box d)$\hspace{2.5em}\=
$\varepsilon^\Box_A\cirk d=\mj_A$,\hspace{7em}\=
$(\varepsilon^\Diamond g)$\hspace{2.5em}\=
$g\cirk\varepsilon^\Diamond_A=\mj_A$, \hspace{1.7em}\=
$(\varepsilon^\Box d)$\hspace{2.5em}\= $\varepsilon^\Box_A\cirk
d=\mj_A$,\hspace{7em}\= $(\varepsilon^\Diamond g)$\hspace{2.5em}\=
$g\cirk\varepsilon^\Diamond_A=\mj_A$,\kill

\>($\delta^{\Box M}$~{\it nat~h})\hspace{1em}$\Box h\cirk
d_1=d_2\cirk h$,\>\>($\delta^{\Diamond M}$~{\it
nat~h})\hspace{1em} $g_2\cirk\Diamond h=h\cirk g_1$.
\end{tabbing}

We define two functors $R$ and $L$ from $\cal C$ to ${\cal
C}^\Box_\Diamond$ in the following manner. The object $RA$ is
${\langle\Box
A,\delta^{\Box\Box}_A,\delta^{\Diamond\Box}_A\rangle}$, while $Rf$
is $\Box f$, appropriately indexed. Dually, $LA$ is
${\langle\Diamond
A,\delta^{\Box\Diamond}_A,\delta^{\Diamond\Diamond}_A\rangle}$,
while $Lf$ is $\Diamond f$, appropriately indexed. We define next
a functor $U$ from ${\cal C}^\Box_\Diamond$ to $\cal C$ by
stipulating that ${U\langle A,d,g\rangle}$ is $A$ and $Uh$ is $h$.
Then it can be shown that $L$ is left adjoint to $U$, while $R$ is
right adjoint to $U$. We need the equation ($\delta${\lat N}$\;g$)
to check that the counit of the adjunction involving $L$ and $U$
satisfies ($\delta^{\Box M}$~{\it nat~h}). Dually, we need the
equation ($\delta${\cir I}$\;d$) to check that the unit of the
adjunction involving $U$ and $R$ satisfies ($\delta^{\Diamond
M}$~{\it nat~h}). The endofunctors $UR$ and $UL$ are equal
respectively to $\Box$ and $\Diamond$.

We have a trijunction with the categories ${\cal C}^\Box_\Diamond$
and $\cal C$ above, together with the functors $L$, $R$ and $U$
between them, and the dyad to which this trijunction gives rise is
the dyad on $\cal C$. One can prove a theorem that says that this
trijunction is terminal, in an appropriate sense, among all the
trijunctions that give rise to the dyad on $\cal C$, which is
analogous to a theorem about the adjunction involving the
Eilenberg-Moore category (see \cite{ML98}, Section VI.3, and
\cite{D99}, Section 5.2.4).

Consider the full subcategory $({\cal C}^\Box_\Diamond)_{\mbox{\it
free}}$ of ${\cal C}^\Box_\Diamond$ whose objects are of the form
$\langle MA,\delta^{\Box M}_A,\delta^{\Diamond M}_A\rangle$. It is
clear that there is a trijunction involving $({\cal
C}^\Box_\Diamond)_{\mbox{\it free}}$ and $\cal  C$, but it is not
immediately clear how to obtain from $({\cal
C}^\Box_\Diamond)_{\mbox{\it free}}$ an analogue of the Kleisli
category, such that the trijunction involving it and $\cal C$
would be initial among all the trijunctions that give rise to the
dyad on $\cal C$ (see \cite{ML98}, Section VI.5, and \cite{D99},
Sections 5.1.6 and 5.2.4). We leave this matter for another
occasion.

We can prove coherence for trijunctions with respect to a functor
$G$ into \emph{Gen} such that the counits and units of the
trijunction are mapped into the split equivalences corresponding
to the following pictures:
\begin{center}
\begin{picture}(240,30)(0,60)
\put(74,78){\line(1,0){5}} \put(74,69){\line(1,0){5}}
\put(74,69){\line(0,1){9}} \put(79,69){\line(0,1){9}}
\put(64,78){\oval(10,10)[b]}

\put(80,80){\makebox(0,0)[br]{$LUA$}}
\put(80,60){\makebox(0,0)[br]{$A$}}
\put(30,71){\makebox(0,0)[r]{$G\varphi^L_A$}}

\put(224,78){\line(1,0){5}} \put(224,69){\line(1,0){5}}
\put(224,69){\line(0,1){9}} \put(229,69){\line(0,1){9}}
\put(214,78){\oval(10,10)[b]}

\put(230,80){\makebox(0,0)[br]{$URB$}}
\put(230,60){\makebox(0,0)[br]{$B$}}
\put(180,71){\makebox(0,0)[r]{$G\varphi^R_B$}}

\end{picture}
\end{center}

\begin{center}
\begin{picture}(240,30)

\put(74,18){\line(1,0){5}} \put(74,9){\line(1,0){5}}
\put(74,9){\line(0,1){9}} \put(79,9){\line(0,1){9}}
\put(64,9){\oval(10,10)[t]}

\put(80,20){\makebox(0,0)[br]{$A$}}
\put(80,0){\makebox(0,0)[br]{$RUA$}}
\put(30,11){\makebox(0,0)[r]{$G\gamma^R_A$}}

\put(224,18){\line(1,0){5}} \put(224,9){\line(1,0){5}}
\put(224,9){\line(0,1){9}} \put(229,9){\line(0,1){9}}
\put(214,9){\oval(10,10)[t]}

\put(230,20){\makebox(0,0)[br]{$B$}}
\put(230,0){\makebox(0,0)[br]{$ULB$}}
\put(180,11){\makebox(0,0)[r]{$G\gamma^L_B$}}

\end{picture}
\end{center}
and, for $F$ being $U$, $L$ or $R$, we have
\begin{center}
\begin{picture}(80,30)
\put(68,17.5){\line(1,0){9}} \put(68,9){\line(1,0){9}}
\put(68,9){\line(0,1){8.5}} \put(77,9){\line(0,1){8.5}}
\put(63,9){\line(0,1){8.5}}

\put(80,20){\makebox(0,0)[br]{$FC_1$}}
\put(80,0){\makebox(0,0)[br]{$FC_2$}}
\put(30,11){\makebox(0,0)[r]{$G(Ff)$}}
\put(77,13){\makebox(0,0)[r]{\tiny$G\!f$}}

\end{picture}
\end{center}
(Related functors may be found in \cite{D99}, Section 4.10,
\cite{D08}, Section~7, \cite{DP03c} and \cite{DP08a}, Section~6;
in contradistinction to what we have in \cite{DP03c} and
\cite{DP08a}, circles cannot arise with trijunctions, as they do
not arise in \cite{D99} and \cite{D08}). The image of this functor
$G$ is included in a subcategory of \emph{Gen} called \emph{Br} in
\cite{DP07} (Section 2.3), where the members of the partitions
induced by the split equivalences are two-element sets. To prove
this coherence result we can rely on a normal form $f_2\cirk f_1$
for the arrow terms of freely generated trijunctions where,
besides $U$, $L$, and $R$, we find in $f_1$ only $\varphi^R$ and
$\varphi^L$, and in $f_2$ only $\gamma^L$ and $\gamma^R$ (see
\cite{D99}, Chapter~4, and \cite{D08}, Sections 5 and 7, for an
analogous result for adjunctions).

Our coherence results for $S5_{\Box\Diamond}$ and
$5S_{\Box\Diamond}$, established in the preceding sections, are
closely related to this coherence result for trijunctions. The
connection of the functors $G$ from $S5_{\Box\Diamond}$ and
$5S_{\Box\Diamond}$ to \emph{Gen} with the functor $G$ for
trijunctions is explained in \cite{DP08a} (end of Section~6). The
trijunctional split equivalences are an isomorphic image of the
split equivalences of $S5_{\Box\Diamond}$ and $5S_{\Box\Diamond}$.
(In the terminology of \cite{DP08a}, Sections 6-7, the split
equivalences of $S5_{\Box\Diamond}$ and $5S_{\Box\Diamond}$ arise
out of the even equivalence classes, i.e.\ the black regions, of
trijunctional split equivalences.)

If we generate freely a trijunction with a single generating
object of the category $\cal B$, then $\cal B$ is isomorphic to
the free dyad generated by a single object, i.e.\ to the category
$S5_{\Box\Diamond}$. If we generate freely our trijunction with a
single generating object of the category $\cal A$, then $\cal A$
is isomorphic to the free codyad generated by a single object,
i.e.\ to the category $5S_{\Box\Diamond}$. This is shown by
relying on the coherence results for trijunctions,
$S5_{\Box\Diamond}$ and $5S_{\Box\Diamond}$ mentioned in the
preceding paragraph. Related matters are considered at the end of
the paper in connection with the square of trijunctions (see
Section 11).

\section{Maximality in the context of $S4$}
Let $S4_{\Box triv}$ be the category defined like $S4_\Box$ save
that for every $n$ we have the additional equation
\begin{tabbing}
\hspace{1.7em}($\varepsilon^\Box$~{\it
triv})\hspace{9.7em}$\Box\varepsilon^\Box_n=\varepsilon^\Box_{n+1}$.
\end{tabbing}
It is shown in \cite{D99} (Section 5.8.2) that the category
$S4_{\Box triv}$ is a \emph{preorder}; namely, for every $f$ and
$g$ of the same type we have ${f=g}$. In $S4_{\Box triv}$ we have
that $\Box$ is isomorphic to $\Box\Box$.

To define $S4_{\Box triv}$, we could use instead of
($\varepsilon^\Box$~{\it triv}) the equation
\[
\Box\delta^{\Box\Box}_n=\delta^{\Box\Box}_{n+1},
\]
which would make superfluous the assumption of the equation
$(\delta^{\Box\Box})$. As a matter of fact, to define $S4_{\Box
triv}$, we could add to $S4_\Box$ instead of
($\varepsilon^\Box$~{\it triv}) any other equation between arrow
terms of $S4_\Box$ that does not hold in $S4_\Box$, provided we
assume this equation \emph{universally}. This means that besides
this equation we assume also all the equations obtained from it by
increasing the subscripts of $\mj$, $\varepsilon^\Box$ and
$\delta^{\Box\Box}$ by a natural number $k$. For example, if we
assume the following instance of ($\varepsilon^\Box$~{\it triv}):
\[
\Box\varepsilon^\Box_1=\varepsilon^\Box_2,
\]
we must assume also
${\Box\varepsilon^\Box_{1+k}=\varepsilon^\Box_{2+k}}$ for every
${k\geq 0}$. We do not assume thereby
${\Box\varepsilon^\Box_0=\varepsilon^\Box_1}$, but it can be shown
that this last instance of ($\varepsilon^\Box$~{\it triv}) is
derivable from $\Box\varepsilon^\Box_1=\varepsilon^\Box_2$, and so
we obtain the whole of ($\varepsilon^\Box$~{\it triv}).

In defining the categories of this paper we always assume
universally the axiomatic equations. So when for the extensions we
assume universally new equations, we proceed as usual in our
definitions.

The \emph{maximality} of $S4_\Box$ is the result which says that
any extension of the definition of that category with a new
universally-holding equation for the arrow terms of that category
(\emph{new} meaning that it does not hold in $S4_\Box$) leads to
collapse, i.e.\ to a category that is a preorder. (For a proof of
this result, see \cite{D99}, Section 5.10.) We will speak of
maximality for other categories later on in the same sense. (The
notion of maximality in \cite{DP04}, Section 9.3, is related, but
stronger; it requires not only that the newly obtained category,
like $S4_{\Box triv}$, be a preorder, but also that any category
in the class in which the newly obtained category is the freely
generated one be also a preorder.)

{\samepage The category $S4_{\Diamond triv}$ is defined like
$S4_\Diamond$ with the additional equation
\begin{tabbing}
\hspace{1.7em}($\varepsilon^\Diamond$~{\it
triv})\hspace{9.5em}$\Diamond\varepsilon^\Diamond_n=\varepsilon^\Diamond_{n+1}$.
\end{tabbing}}
\noindent We can say for $S4_{\Diamond triv}$, \emph{mutatis
mutandis}, whatever we said for $S4_{\Box triv}$. The category
$S4_\Diamond$ is maximal in the same sense in which $S4_\Box$ is
maximal.

When we consider extensions with new equations for categories like
$S4_{\Box\Diamond}$, whose objects are not finite ordinals but
modalities, assuming an equation \emph{universally} means that
besides this equation we assume also all the equations obtained
from it by appending to the subscripts of the primitive arrow
terms an arbitrary modality $A$ on the right-hand side. For
example, the equation ($\varepsilon^\Box$~{\it triv}) now becomes
the following equation:
\begin{tabbing}
\hspace{1.7em}($\varepsilon^\Box$~{\it
triv})\hspace{9.5em}$\Box\varepsilon^\Box_A=\varepsilon^\Box_{\Box
A}$.
\end{tabbing}
If we assume the following instance of this equation:
\[
\Box\varepsilon^\Box_\Diamond=\varepsilon^\Box_{\Box\Diamond},
\]
we must assume also ${\Box\varepsilon^\Box_{\Diamond
A}=\varepsilon^\Box_{\Box\Diamond A}}$ for every modality $A$.

The category $S4_{\Box\Diamond}$ is not maximal in the sense in
which $S4_\Box$ and $S4_\Diamond$ were maximal. We can add to
$S4_{\Box\Diamond}$ one of the equations ($\varepsilon^\Box$~{\it
triv}) or ($\varepsilon^\Diamond$~{\it triv}), where $n$ and
${n\pl 1}$ are replaced respectively by $A$ and $\Box A$, or $A$
and $\Diamond A$, without thereby obtaining the other. This is
shown with the help of appropriate modifications of the functor
$G$ from $S4_{\Box\Diamond}$ to \emph{Rel} (we may omit the pairs
involving $\Box$ without omitting those involving $\Diamond$, and
vice versa).

Let the category $S4_{\Box\Diamond\sharp}$ be defined like
$S4_{\Box\Diamond}$ save that we have the additional equations
($\varepsilon^\Box$~{\it triv}) and ($\varepsilon^\Diamond$~{\it
triv}), with the replacement mentioned in the preceding paragraph.
This category is not a preorder because the equation
\begin{tabbing}
\hspace{1.7em}
$(\Box\Diamond)$\hspace{6.7em}$\Diamond\Box\varepsilon^\Diamond_A\cirk\varepsilon^\Box_{\Diamond\Box
A}=\varepsilon^\Diamond_{\Box\Diamond
A}\cirk\Box\Diamond\varepsilon^\Box_A$,
\end{tabbing}
does not hold in it, as we are going to show now. Consider the
pictures
\begin{center}
\begin{picture}(240,40)

\put(61,9){\line(2,5){7}} \put(68,9){\line(2,5){7}}

\put(80,29){\makebox(0,0)[br]{$\Box\Diamond\Box$}}
\put(80,0){\makebox(0,0)[br]{$\Diamond\Box\Diamond$}}
\put(35,16){\makebox(0,0)[r]{$G(\Diamond\Box\varepsilon^\Diamond\cirk\varepsilon^\Box_{\Diamond\Box})$}}

\put(226,9){\line(-2,5){7}} \put(218,9){\line(-2,5){7}}

\put(230,29){\makebox(0,0)[br]{$\Box\Diamond\Box$}}
\put(230,0){\makebox(0,0)[br]{$\Diamond\Box\Diamond$}}
\put(185,16){\makebox(0,0)[r]{$G(\varepsilon^\Diamond_{\Box\Diamond}\cirk\Box\Diamond\varepsilon^\Box)$}}

\end{picture}
\end{center}
which are yielded by the functor $G$ from $S4_{\Box\Diamond}$ to
\emph{Rel}, but also by a modification $G^\sharp$ of that functor,
which goes from $S4_{\Box\Diamond\sharp}$ to \emph{Rel}, and takes
into account that $MM$ is isomorphic to $M$, for $M$ being $\Box$
or $\Diamond$.

To define $G^\sharp$, we define first inductively a function
$^\sharp$ on the objects of $S4_{\Box\Diamond\sharp}$, which are
also objects of $S4_{\Box\Diamond}$, i.e.\ the modalities. For
$M,M_1,M_2\in\{\Box,\Diamond\}$ we have
\[
\begin{array}{c}
M^\sharp=M,
\\[2ex]
(M_1M_2A)^\sharp=\left\{ \begin{array}{ll} (M_2A)^\sharp &
\mbox{if $M_1$ is $M_2$,}
\\
M_1(M_2A)^\sharp & \mbox{if $M_1$ is not $M_2$.}
\end{array}
\right.
\end{array}
\]
Next we define inductively the arrow terms ${j_A\!:A\vdash
A^\sharp}$ and ${j^A\!:A^\sharp\vdash A}$ of $S4_{\Box\Diamond}$:
\begin{tabbing}
\centerline{$j_M=j^M=\mj_M$,}
\\[1.5ex]
\hspace{7em}\= $j_{\Box\Box A}=j_{\Box
A}\cirk\varepsilon^\Box_{\Box A}$,\hspace{5em}\= $j^{\Box\Box
A}=\delta^{\Box\Box}_A\cirk j^{\Box A}$,
\\*[1ex]
\> $j_{\Diamond\Diamond A}=j_{\Diamond
A}\cirk\delta^{\Diamond\Diamond}_A$,\> $j^{\Diamond\Diamond
A}=\varepsilon^\Diamond_{\Diamond A}\cirk j^{\Diamond A}$,
\\[1.5ex]
\hspace{1em}for $M_1$ different from $M_2$,
\\*
\> $j_{M_1M_2A}=M_1j_{M_2A}$,\> $j^{M_1M_2A}=M_1j^{M_2A}$.
\end{tabbing}
It is easy to see that $j_A$ and $j^A$ are isomorphisms of
$S4_{\Box\Diamond\sharp}$, inverse to each other. Then, for $G$
being the functor from $S4_{\Box\Diamond}$ to \emph{Rel}, we have
that $G^\sharp A$ is $GA^\sharp$, and for ${f\!:A\vdash B}$ an
arrow term of $S4_{\Box\Diamond\sharp}$, i.e.\ of
$S4_{\Box\Diamond}$, we have that $G^\sharp f$ is ${G(j_B\cirk
f\cirk j^A)}$. It is easy to verify that $G^\sharp$ is indeed a
functor, which is sufficient to show that the equation
$(\Box\Diamond)$ does not hold in $S4_{\Box\Diamond\sharp}$.

Then we can infer from $S4_{\Box\Diamond}$ Coherence that
$G^\sharp$ is a faithful functor, i.e.\ $S4_{\Box\Diamond\sharp}$
\emph{Coherence}. Suppose for ${f,g\!:A\vdash B}$ that ${G^\sharp
f=G^\sharp g}$; by $S4_{\Box\Diamond}$ Coherence we have $j_B\cirk
f\cirk j^A=j_B\cirk g\cirk j^A$ in $S4_{\Box\Diamond}$, and hence
also in $S4_{\Box\Diamond\sharp}$. Since $j^A$ and $j_B$ are
isomorphisms in $S4_{\Box\Diamond\sharp}$, it follows that ${f=g}$
in $S4_{\Box\Diamond\sharp}$.

Let $S4_{\Box\Diamond triv}$ be defined like
$S4_{\Box\Diamond\sharp}$ save that we have the additional
equation $(\Box\Diamond)$. In $S4_{\Box\Diamond triv}$, besides
having that $MM$ is isomorphic to $M$, for $M$ being $\Box$ or
$\Diamond$, we also have this isomorphism for $M$ being
$\Box\Diamond$ or $\Diamond\Box$. For $M$ being $\Box\Diamond$,
let
\begin{tabbing}
\hspace{9em}\= $\;\;\;\;i$ \=
$=_{df}\Box\delta^{\Diamond\Diamond}\cirk\Box\Diamond\varepsilon^\Box_\Diamond\!:$
\= $MM\vdash M$,
\\*[1.5ex]
\> $i^{-1}$\>
$=_{df}\Box\varepsilon^\Diamond_{\Box\Diamond}\cirk\delta^{\Box\Box}_\Diamond\!:$
\> $M\vdash MM$.
\end{tabbing}
To show that ${i\cirk i^{-1}=\mj_{\Box\Diamond}}$, we may apply
$S4_{\Box\Diamond}$ Coherence. To show that ${i^{-1}\cirk
i=\mj_{\Box\Diamond\Box\Diamond}}$, we have
\begin{tabbing}
\hspace{3em} $i^{-1}\cirk i$ \=
$=\Box\Diamond\Box\delta^{\Diamond\Diamond}\cirk\Box\varepsilon^\Diamond_{\Box\Diamond\Diamond}
\cirk\Box\Box\Diamond\varepsilon^\Box_\Diamond\cirk\delta^{\Box\Box}_{\Diamond\Box\Diamond}$,
\hspace{.5em}by $S4_{\Box\Diamond}$ Coherence,
\\*[1.5ex]
\>
$=\Box\Diamond\Box\delta^{\Diamond\Diamond}\cirk\Box\Diamond\Box\varepsilon^\Diamond_\Diamond\cirk
\Box\varepsilon^\Box_{\Diamond\Box\Diamond}\cirk\delta^{\Box\Box}_{\Diamond\Box\Diamond}$,
\hspace{.5em}by $(\Box\Diamond)$,
\\[1.5ex]
\> $=\mj_{\Box\Diamond\Box\Diamond}$, \hspace{.5em}by
$S4_{\Box\Diamond}$ Coherence.
\end{tabbing}
We proceed analogously for $M$ being $\Diamond\Box$.

For the proposition below, we need the following diagram of arrows
of $S4_{\Box\Diamond triv}$, which without the arrow terms may be
found in \cite{K22}, and is commonly used to classify the
modalities of $S4$ (see \cite{HC96}, p.\ 56):
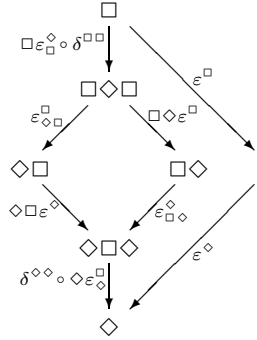
\begin{center}
\begin{picture}(90,137)

\put(30,30){\vector(0,-1){17}} \put(30,120){\vector(0,-1){17}}
\put(22,90){\vector(-1,-1){17}} \put(38,90){\vector(1,-1){17}}
\put(55,60){\vector(-1,-1){17}} \put(5,60){\vector(1,-1){17}}
\put(85,60){\vector(-1,-1){47}} \put(38,120){\vector(1,-1){47}}

\put(30,5){\makebox(0,0){$\Diamond$}}
\put(30,35){\makebox(0,0){$\Diamond\Box\Diamond$}}
\put(0,65){\makebox(0,0){$\Diamond\Box$}}
\put(60,65){\makebox(0,0){$\Box\Diamond$}}
\put(30,95){\makebox(0,0){$\Box\Diamond\Box$}}
\put(30,125){\makebox(0,0){$\Box$}}

\put(2,48){\makebox(0,0)[b]{\scriptsize$\Diamond\Box\varepsilon^\Diamond$}}
\put(54,48){\makebox(0,0)[b]{\scriptsize$\varepsilon^\Diamond_{\Box\Diamond}$}}

\put(7,84){\makebox(0,0)[b]{\scriptsize$\varepsilon^\Box_{\Diamond\Box}$}}
\put(54.5,84){\makebox(0,0)[b]{\scriptsize$\Box\Diamond\varepsilon^\Box$}}

\put(13,111){\makebox(0,0)[b]{\scriptsize$\Box\varepsilon^\Diamond_\Box\!\cirk\delta^{\Box\Box}$}}
\put(13,22){\makebox(0,0)[b]{\scriptsize$\delta^{\Diamond\Diamond}\!\cirk
\Diamond\varepsilon^\Box_\Diamond$}}

\put(66,98){\makebox(0,0)[b]{\scriptsize$\varepsilon^\Box$}}
\put(66,31.5){\makebox(0,0)[b]{\scriptsize$\varepsilon^\Diamond$}}

\end{picture}
\end{center}

\prop{Preorder of $S4_{\Box\Diamond triv}$}{The category
$S4_{\Box\Diamond triv}$ is a preorder, and its skeleton is given
by the diagram above.}

\dkz Note first that the isomorphisms of $S4_{\Box\Diamond triv}$
yield just the seven objects in the diagram above. Next, for
${(M_1,M_2)}$ being a pair of these seven modalities that is not
${(\Box\Diamond\Box,\Diamond\Box\Diamond)}$, we may conclude from
$S4_{\Box\Diamond}$ Coherence, and the properties of the functor
$G$ from $S4_{\Box\Diamond}$ to \emph{Rel}, that there is at most
one arrow from $M_1$ to $M_2$ in $S4_{\Box\Diamond}$, and hence
also in $S4_{\Box\Diamond triv}$. (Every occurrence of $\Box$ in
the target is linked to an occurrence of $\Box$ in the source, and
every occurrence of $\Diamond$ in the source is linked to an
occurrence of $\Diamond$ in the target; moreover, links are not
crossed with each other.) There are two arrows from
$\Box\Diamond\Box$ to $\Diamond\Box\Diamond$ in
$S4_{\Box\Diamond}$, which make the two paths in the small square
in the diagram above. They are instances of the two sides of
$(\Box\Diamond)$. So all paths in the diagram above commute in
$S4_{\Box\Diamond triv}$.\qed

\vspace{2ex}

We can then show the following.

\prop{Maximality of $S4_{\Box\Diamond\sharp}$}{The category
$S4_{\Box\Diamond\sharp}$ is maximal.}

\dkz Suppose we have the arrow terms ${f,g\!:A\vdash B}$ of
$S4_{\Box\Diamond}$ such that ${f=g}$ does not hold in
$S4_{\Box\Diamond\sharp}$. By $S4_{\Box\Diamond\sharp}$ Coherence,
we have $G^\sharp f\neq G^\sharp g$, where $G^\sharp$ is the
functor from $S4_{\Box\Diamond\sharp}$ to \emph{Rel} defined
above. Then it can be inferred that $G^\sharp f$ corresponds to
the picture with solid lines, while $G^\sharp g$ corresponds to
the picture with dotted lines
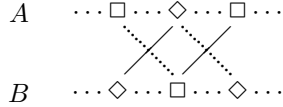
\begin{center}
\begin{picture}(120,45)

\put(40,12){\line(1,1){18}} \put(62,12){\line(1,1){18}}

\multiput(58,12)(-1.5,1.5){13}{\circle*{.3}}
\multiput(80,12)(-1.5,1.5){13}{\circle*{.3}}

\put(60,5){\makebox(0,0){$\cdots\Diamond\cdots\Box\cdots\Diamond\cdots$}}
\put(60,35){\makebox(0,0){$\cdots\Box\cdots\Diamond\cdots\Box\cdots$}}
\put(0,5){\makebox(0,0){$B$}} \put(0,35){\makebox(0,0){$A$}}
\end{picture}
\end{center}
($G^\sharp f$ and $G^\sharp g$ can of course switch places). This
is because in our pictures we cannot have crossings. Let
${h_A\!:\Box\Diamond\Box\vdash A}$ and
${h_B\!:B\vdash\Box\Diamond\Box}$ be the arrows of
$S4_{\Box\Diamond\sharp}$ such that $G^\sharp h_A$ and $G^\sharp
h_B$ correspond respectively to the pictures
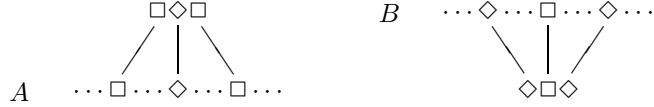
\begin{center}
\begin{picture}(240,45)

\put(60,12){\line(0,1){18}} \put(39,12){\line(2,3){12}}
\put(81,12){\line(-2,3){12}}

\put(200,12){\line(0,1){18}} \put(179,30){\line(2,-3){12}}
\put(221,30){\line(-2,-3){12}}

\put(200,35){\makebox(0,0){$\cdots\Diamond\cdots\Box\cdots\Diamond\cdots$}}
\put(60,5){\makebox(0,0){$\cdots\Box\cdots\Diamond\cdots\Box\cdots$}}
\put(140,35){\makebox(0,0){$B$}} \put(0,5){\makebox(0,0){$A$}}
\put(60,35){\makebox(0,0){$\Box\Diamond\Box$}}
\put(200,5){\makebox(0,0){$\Diamond\Box\Diamond$}}

\end{picture}
\end{center}
In the left picture, any $\Box$ in $A$ to the left of the
displayed $\Diamond$ is tied to the left $\Box$ in
$\Box\Diamond\Box$, and analogously when ``left'' is replaced by
``right''. We interpret the right picture analogously replacing
$\Box$ by $\Diamond$. Then, by $S4_{\Box\Diamond\sharp}$
Coherence, we can conclude that
\begin{tabbing}
\hspace{11.5em}\= $h_B\cirk f\cirk h_A\,$\=
$=\Diamond\Box\varepsilon^\Diamond\cirk\varepsilon^\Box_{\Diamond\Box}$,
\\*[1.5ex]
\> $h_B\cirk g\cirk h_A$ \>
$=\varepsilon^\Diamond_{\Box\Diamond}\cirk\Box\Diamond\varepsilon^\Box$,
\end{tabbing}
and this, together with appending modalities on the right-hand
side in the subscripts of $f$ and $g$, yields the equation
$(\Box\Diamond)$. So if we assume ${f=g}$ universally, we will
also have $(\Box\Diamond)$, and hence we will be in
$S4_{\Box\Diamond triv}$, which is a preorder.
\mbox{\hspace{1em}}\qed

\vspace{2ex}

Let the category $S4.2_{\Box\Diamond\sharp}$ be defined like
$S4.2_{\Box\Diamond}$ save that we have the additional equations
($\varepsilon^\Box$~{\it triv}) and ($\varepsilon^\Diamond$~{\it
triv}) we used to obtain $S4_{\Box\Diamond\sharp}$ out of
$S4_{\Box\Diamond}$. We define a functor $G^\sharp$ from
$S4.2_{\Box\Diamond\sharp}$ to \emph{Rel} with the help of the
functor $G$ from $S4.2_{\Box\Diamond}$ to \emph{Rel}, as we did
for $S4_{\Box\Diamond\sharp}$, and by relying on
$S4.2_{\Box\Diamond}$ Coherence we establish that this new functor
$G^\sharp$ is faithful, i.e.\ $S4.2_{\Box\Diamond\sharp}$
\emph{Coherence} (see above). Then we can show that
$S4.2_{\Box\Diamond\sharp}$ is not a preorder. Besides the
equation $(\Box\Diamond)$, we do not have in
$S4.2_{\Box\Diamond\sharp}$ the equations
\begin{tabbing}
\hspace{1em}\= $(\Box\varepsilon\chi)$\hspace{1em}\=
$\Box\Diamond\varepsilon^\Box_A=\chi^{\Diamond\Box}_A\cirk\varepsilon^\Box_{\Diamond\Box
A}$,\hspace{5em}\= $(\Diamond\varepsilon\chi)$\hspace{1em}\=
$\Diamond\Box\varepsilon^\Diamond_A=\varepsilon^\Diamond_{\Box\Diamond
A}\cirk\chi^{\Diamond\Box}_A$,
\\[1.5ex]
\> $(\Box\varepsilon\chi\delta)$\> $\chi^{\Diamond\Box}_{\Box
A}\cirk
\Diamond\delta^{\Box\Box}_A\cirk\varepsilon^\Box_{\Diamond\Box
A}=\mj_{\Box\Diamond\Box A}$,\>
$(\Diamond\varepsilon\chi\delta)$\>
$\varepsilon^\Diamond_{\Box\Diamond
A}\cirk\Box\delta^{\Diamond\Diamond}_A\cirk\chi^{\Diamond\Box}_{\Diamond
A}=\mj_{\Diamond\Box\Diamond A}$,
\end{tabbing}
as it is easily shown with the help of $G^\sharp$.

When we add $(\Box\Diamond)$ to $S4.2_{\Box\Diamond}$, we can
derive $(\Box\varepsilon\chi)$ as follows:
\begin{tabbing}
\hspace{5em}$\Box\Diamond\varepsilon^\Box_A$ \=
$=\Box\delta^{\Diamond\Diamond}_A\cirk\Box\varepsilon^\Diamond_{\Diamond
A}\cirk\Box\Diamond\varepsilon^\Box_A$, \hspace{.5em}by
$(\Diamond\Diamond\beta)$,
\\[1.5ex]
\>
$=\Box\delta^{\Diamond\Diamond}_A\cirk\chi^{\Diamond\Box}_{\Diamond
A}\cirk\varepsilon^\Diamond_{\Box\Diamond A}\cirk
\Box\Diamond\varepsilon^\Box_A$, \hspace{.5em}by
$(\varepsilon^\Diamond\chi^{\Diamond\Box})$,
\\[1.5ex]
\>
$=\Box\delta^{\Diamond\Diamond}_A\cirk\chi^{\Diamond\Box}_{\Diamond
A}\cirk \Diamond\Box\varepsilon^\Diamond_A\cirk
\varepsilon^\Box_{\Diamond\Box A}$, \hspace{.5em}by
$(\Box\Diamond)$,
\\[1.5ex]
\>
$=\Box\delta^{\Diamond\Diamond}_A\cirk\Box\Diamond\varepsilon^\Diamond_A\cirk\chi^{\Diamond\Box}_A\cirk
\varepsilon^\Box_{\Diamond\Box A}$, \hspace{.5em}by
($\chi^{\Diamond\Box}$~{\it nat}),
\\[1.5ex]
\> $=\chi^{\Diamond\Box}_A\cirk \varepsilon^\Box_{\Diamond\Box
A}$, \hspace{.5em}by $(\Diamond\Diamond\eta)$.
\end{tabbing}
(As a matter of fact, by $S4.2_{\Box\Diamond}$ Coherence we can
pass immediately to the second line, and also from the third line
to the last line.) We proceed analogously to derive
$(\Diamond\varepsilon\chi)$ from $(\Box\Diamond)$.

Next, when we add $(\Box\varepsilon\chi)$ to
$S4.2_{\Box\Diamond}$, we can derive $(\Box\varepsilon\chi\delta)$
as follows:
\begin{tabbing}
\hspace{5em}$\chi^{\Diamond\Box}_{\Box
A}\cirk\Diamond\delta^{\Box\Box}_A\cirk\varepsilon^\Box_{\Diamond\Box
A}$ \= $=\chi^{\Diamond\Box}_{\Box
A}\cirk\varepsilon^\Box_{\Diamond\Box\Box
A}\cirk\Box\Diamond\delta^{\Box\Box}_A$, \hspace{.5em}by
($\varepsilon^\Box$~{\it nat}),
\\[1.5ex]
\> $=\Box\Diamond\varepsilon^\Box_{\Box
A}\cirk\Box\Diamond\delta^{\Box\Box}_A$, \hspace{.5em}by
$(\Box\varepsilon\chi)$,
\\[1.5ex]
\> $=\mj_{\Box\Diamond\Box A}$, \hspace{.5em}by $(\Box\Box\beta)$.
\end{tabbing}
We proceed analogously to derive $(\Diamond\varepsilon\chi\delta)$
from $(\Diamond\varepsilon\chi)$.

When we add $(\Box\varepsilon\chi\delta)$ and
$(\Diamond\varepsilon\chi\delta)$ to $S4.2_{\Box\Diamond}$, we can
derive $(\Box\Diamond)$ as follows:
\begin{tabbing}
\hspace{1.6em}$\Diamond\Box\varepsilon^\Diamond_A\cirk\varepsilon^\Box_{\Diamond\Box
A}$ \= $=\varepsilon^\Diamond_{\Box\Diamond A}\cirk
\Box\delta^{\Diamond\Diamond}_A\cirk\chi^{\Diamond\Box}_{\Diamond
A}\cirk
\Diamond\Box\varepsilon^\Diamond_A\cirk\varepsilon^\Box_{\Diamond\Box
A}$, \hspace{.5em}by $(\Diamond\varepsilon\chi\delta)$,
\\[1.5ex]
\> $=\varepsilon^\Diamond_{\Box\Diamond A}\cirk
\chi^{\Diamond\Box}_A\cirk \varepsilon^\Box_{\Diamond\Box A}$,
\hspace{.5em}by ($\chi^{\Diamond\Box}$~{\it nat}) and
$(\Diamond\Diamond\eta)$,
\\[1.5ex]
\> $=\varepsilon^\Diamond_{\Box\Diamond A}\cirk
\Box\Diamond\varepsilon^\Box_A\cirk \chi^{\Diamond\Box}_{\Box
A}\cirk \Diamond\delta^{\Box\Box}_A\cirk
\varepsilon^\Box_{\Diamond\Box A}$, \hspace{.5em}by
$(\Box\Box\eta)$ and ($\chi^{\Diamond\Box}$~{\it nat}),
\\[1.5ex]
\> $=\varepsilon^\Diamond_{\Box\Diamond A}\cirk
\Box\Diamond\varepsilon^\Box_A$, \hspace{.5em}by
$(\Box\varepsilon\chi\delta)$.
\end{tabbing}
With the help of modifications of $G^\sharp$ in which we omit all
$\Box$-\emph{links}, i.e.\ links involving $\Box$ (which are here
links joining occurrences of $\Box$), without omitting
$\Diamond$-\emph{links}, i.e.\ links involving $\Diamond$ (which
are here links joining occurrences of $\Diamond$), and vice versa,
we can show that none of $(\Box\varepsilon\chi)$ and
$(\Diamond\varepsilon\chi)$ implies the other, and the same for
$(\Box\varepsilon\chi\delta)$ and
$(\Diamond\varepsilon\chi\delta)$.

Let the category $S4.2_{\Box\Diamond triv}$ be defined like
$S4.2_{\Box\Diamond\sharp}$ save that we have the additional
equation $(\Box\Diamond)$. We can show that $S4.2_{\Box\Diamond
triv}$ is a preorder, and that its skeleton is given by the
following diagram:
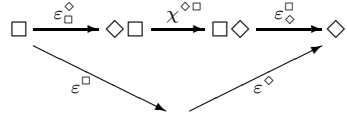
\begin{center}
\begin{picture}(120,45)

\put(5.5,36){\vector(1,0){24.5}} \put(90,36){\vector(1,0){24.5}}
\put(50,36){\vector(1,0){20}}

\put(5.5,30){\vector(2,-1){50}} \put(64.5,5){\vector(2,1){50}}

\put(0,35){\makebox(0,0){$\Box$}}
\put(40,35){\makebox(0,0){$\Diamond\Box$}}
\put(80,35){\makebox(0,0){$\Box\Diamond$}}
\put(120,35){\makebox(0,0){$\Diamond$}}

\put(17.5,40){\makebox(0,0)[b]{\scriptsize$\varepsilon^\Diamond_\Box$}}
\put(63,40){\makebox(0,0)[b]{\scriptsize$\chi^{\Diamond\Box}$}}
\put(101,40){\makebox(0,0)[b]{\scriptsize$\varepsilon^\Box_\Diamond$}}

\put(24,12){\makebox(0,0)[b]{\scriptsize$\varepsilon^\Box$}}
\put(93,12){\makebox(0,0)[b]{\scriptsize$\varepsilon^\Diamond$}}

\end{picture}
\end{center}

Note that in $S4.2_{\Box\Diamond triv}$ the modalities
$\Box\Diamond\Box$ and $\Diamond\Box$ on the one hand, and
$\Diamond\Box\Diamond$ and $\Box\Diamond$ on the other hand, are
isomorphic. Note also that the arrows $\chi^{\Diamond\Box}_{MA}$
are isomorphisms in $S4.2_{\Box\Diamond triv}$.

It can be shown in extending $S4.2_{\Box\Diamond\sharp}$ that if
we have assumed universally any new equation for arrow terms of
$S4.2_{\Box\Diamond}$, then we will obtain one of the equations
$(\Box\varepsilon\chi)$ and $(\Diamond\varepsilon\chi)$, and hence
also one of the equations $(\Box\varepsilon\chi\delta)$ and
$(\Diamond\varepsilon\chi\delta)$. This is not maximality as we
had for $S4_{\Box\Diamond\sharp}$, but it is not very far from it.
A more precise result, which yields this relative maximality, is
stated as follows.

If the new equation ${f=g}$, which does not hold in
$S4.2_{\Box\Diamond\sharp}$, is such that $G^\sharp f$ differs
from $G^\sharp g$ in the $M$-links, for $M$ being $\Box$ or
$\Diamond$, then we can derive $(M\varepsilon\chi)$ and
$(M\varepsilon\chi\delta)$. If $M$ is $\Box$, then we proceed in a
manner analogous to what we had in the proof of the Maximality of
$S4_{\Box\Diamond\sharp}$, with ${h_A\!:\Box\Diamond\Box\vdash A}$
as there and ${h_B\!:B\vdash\Diamond\Box\Diamond}$ replaced by an
arrow of the type $B\vdash\Box\Diamond$, which is either
$\Box\delta^{\Diamond\Diamond}\cirk\chi^{\Diamond\Box}_\Diamond\cirk
h_B$, or constructed more simply than $h_B$. If $M$ is $\Diamond$,
then we proceed dually by replacing $h_A$.

If $M$ is here only $\Box$, then we cannot derive
$(\Diamond\varepsilon\chi)$ and $(\Diamond\varepsilon\chi\delta)$,
and if it is only $\Diamond$, then we cannot derive
$(\Box\varepsilon\chi)$ and $(\Box\varepsilon\chi\delta)$. If $M$
stands here for both $\Box$ and $\Diamond$, i.e., $G^\sharp f$
differs from $G^\sharp g$ both in $\Box$-links and
$\Diamond$-links, then we can derive $(\Box\Diamond)$. (The point
in the proof of the Maximality of $S4_{\Box\Diamond\sharp}$ is
that $G^\sharp f$ and $G^\sharp g$ cannot differ in $\Box$-links
without differing also in $\Diamond$-links, and vice versa.)

\section{The square of adjunctions}
In this section we consider some elementary facts concerning
adjunctions, which we need for the exposition later on.

That a functor $F$ from $\cal  B$ to $\cal  A$ is \emph{left
adjoint} to a functor $G$ from $\cal  A$ to $\cal  B$
(alternatively, $G$ is \emph{right adjoint} to $F$) means that,
for $I_{\cal B}$ and $I_{\cal A}$ being respectively the identity
functors of $\cal  B$ and $\cal  A$, we have a natural
transformation ${\gamma\!:I_{\cal B}\strt GF}$, the \emph{unit} of
the adjunction, and a natural transformation ${\varphi\!:FG\strt
I_{\cal A}}$, the \emph{counit} of the adjunction, which satisfy
the following \emph{triangular equations} for every object $B$ of
$\cal  B$ and every object $A$ of $\cal  A$:
\[
\varphi_{FB}\cirk
F\gamma_B=\mj_{FB},\hspace{5em}G\varphi_A\cirk\gamma_{GA}=\mj_{GA}.
\]
An \emph{adjunction} is a structure made of such functors $F$ and
$G$, and such natural transformations $\gamma$ and $\varphi$ (for
more details, see \cite{ML98}, Chapter~IV, and \cite{D99},
Chapter~4).

Every adjunction generates four adjunctions involving functor
categories, which we display in the following picture, where left
adjoints have solid arrows, and right adjoints have dotted arrows:
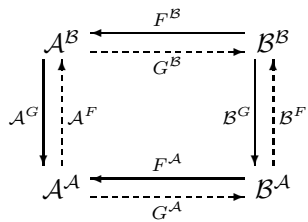
\begin{center}
\begin{picture}(80,80)(0,-10)

\put(-7,53){\vector(0,-1){40}}
\multiput(0,13)(0,4){10}{\line(0,1){2}}
\put(0,51){\vector(0,1){2}}

\put(73,53){\vector(0,-1){40}}
\multiput(80,13)(0,4){10}{\line(0,1){2}}
\put(80,51){\vector(0,1){2}}

\put(69,63){\vector(-1,0){58}}
\multiput(11,56)(4,0){14}{\line(1,0){2}}
\put(67,56){\vector(1,0){2}}

\put(69,8){\vector(-1,0){58}}
\multiput(11,1)(4,0){14}{\line(1,0){2}}
\put(67,1){\vector(1,0){2}}

\put(0,5){\makebox(0,0){${\cal A}^{\cal A}$}}
\put(80,5){\makebox(0,0){${\cal B}^{\cal A}$}}
\put(0,60){\makebox(0,0){${\cal A}^{\cal B}$}}
\put(80,60){\makebox(0,0){${\cal B}^{\cal B}$}}

\put(-8,32){\makebox(0,0)[r]{\scriptsize${\cal A}^G$}}
\put(2,32){\makebox(0,0)[l]{\scriptsize${\cal A}^F$}}
\put(72,32){\makebox(0,0)[r]{\scriptsize${\cal B}^G$}}
\put(82,32){\makebox(0,0)[l]{\scriptsize${\cal B}^F$}}

\put(40,-1){\makebox(0,0)[t]{\scriptsize$G^{\cal A}$}}
\put(40,11){\makebox(0,0)[b]{\scriptsize$F^{\cal A}$}}
\put(40,54){\makebox(0,0)[t]{\scriptsize$G^{\cal B}$}}
\put(40,66){\makebox(0,0)[b]{\scriptsize$F^{\cal B}$}}

\end{picture}
\end{center}

For the functors $H$,$H_1$ and $H_2$ from $\cal  B$ to $\cal  A$,
and for $\alpha$ a natural transformation from $H_1$ to $H_2$, we
have
\begin{tabbing}
\hspace{9em}\= $G^{\cal B}H$\= $=GH$,\hspace{5em}\= $(G^{\cal
B}\alpha)_B$\= $=G\alpha_B$,
\\*[1.5ex]
\> ${\cal A}^GH$\> $=HG$,\> $({\cal A}^G\alpha)_A$\>
$=\alpha_{GA}$;
\end{tabbing}
we define analogously the other functors involved in the
adjunctions above.

In this square of adjunctions, the members of the units for the
two \emph{horizontal} adjunctions are the natural transformations
${\gamma_H\!:H\strt GFH}$, and the members of the counits are
${\varphi_H\!: FGH\strt H}$. For the two \emph{vertical}
adjunctions, the members of the units are ${H\gamma\!:H\strt
HGF}$, and the members of the counits are ${H\varphi\!:HFG\strt
H}$. In the horizontal adjunctions, the functors involving $F$ and
$G$ behave like $F$ and $G$, while in the vertical adjunctions,
the functor involving $F$ becomes right adjoint, and that
involving $G$ left adjoint. The horizontal adjunctions are images
of the original adjunction by two covariant 2-endofunctors of the
2-category \emph{Cat} of categories with functors and natural
transformations, while the vertical adjunctions are such images by
two contravariant 2-endofunctors (for the notions of 2-category
and 2-functor, see \cite{ML98}, Sections XII.3-4).

For ${\cal C}_1,{\cal C}_2\in\{{\cal A},{\cal B}\}$, let a
\emph{canonical} functor from ${\cal C}_1$ to ${\cal C}_2$ be any
functor from ${\cal C}_1$ to ${\cal C}_2$ defined in terms of the
identity functors $I_{\cal A}$ and $I_{\cal B}$, the functors $F$
and $G$, and composition of functors. Let $C{\cal C}_2^{{\cal
C}_1}$ be the subcategory of the functor category ${\cal
C}_2^{{\cal C}_1}$ whose objects are the canonical functors from
${\cal C}_1$ to ${\cal C}_2$, and whose arrows are the
\emph{canonical} natural transformations, defined in terms of the
identity natural transformations, the unit $\gamma$ and counit
$\varphi$ of the adjunction, the functors $F$ and $G$, and
composition. So the objects of $C{\cal B}^{\cal B}$ are $I_{\cal
B}$, $GF$, $GFGF$, etc., those of $C{\cal A}^{\cal B}$ are $F$,
$FGF$, $FGFGF$, etc., those of $C{\cal A}^{\cal A}$ are $I_{\cal
A}$, $FG$, $FGFG$, etc., and finally those of $C{\cal B}^{\cal A}$
are $G$, $GFG$, $GFGF$, etc. Then from the square of adjunctions
above we obtain an analogous square by replacing ${\cal
C}_2^{{\cal C}_1}$ with $C{\cal C}_2^{{\cal C}_1}$. Yet another
analogous square of adjunctions is obtained when $C{\cal
C}_2^{{\cal C}_1}$ is understood as the full subcategory of ${\cal
C}_2^{{\cal C}_1}$ whose objects are the canonical functors from
${\cal C}_1$ to ${\cal C}_2$. (The four preordering equations of
\cite{D99}, Section 4.6.2, are connected by the bijections between
hom-sets of the horizontal and vertical adjunctions in the square
of adjunctions.)

For every category $\cal A$ treated in this paper, whose objects
are either finite ordinals or modalities, let a \emph{canonical}
functor from $\cal A$ to $\cal A$ be a functor definable in terms
of the functors assumed for defining $\cal A$ and composition of
functors. Then these canonical functors may be identified with the
objects of $\cal A$, and, for $C{\cal A}^{\cal A}$ being the full
subcategory of ${\cal A}^{\cal A}$ whose objects are the canonical
functors from $\cal A$ to $\cal A$, we have that $\cal A$ is
isomorphic to $C{\cal A}^{\cal A}$.

If $\cal A$ is $S5_{\Box\Diamond}$, then, as we have seen in
Section~6, the endofunctor $\Diamond$ is left adjoint to the
endofunctor $\Box$. Since $C{\cal A}^{\cal A}$ is isomorphic to
$\cal A$, the $C{\cal C}_2^{{\cal C}_1}$ variant of the square of
adjunctions reduces to
\begin{center}
\begin{picture}(80,75)

\put(-10,53){\vector(0,-1){40}}
\multiput(-3,13)(0,4){10}{\line(0,1){2}}
\put(-3,51){\vector(0,1){2}}

\put(66,64){\vector(-1,0){52}}
\multiput(14,57)(4,0){12}{\line(1,0){2}}
\put(64,57){\vector(1,0){2}}

\put(0,5){\makebox(0,0){$S5_{\Box\Diamond}$}}
\put(0,60){\makebox(0,0){$S5_{\Box\Diamond}$}}
\put(80,60){\makebox(0,0){$S5_{\Box\Diamond}$}}

\put(-11,32){\makebox(0,0)[r]{\scriptsize$I^\Box$}}
\put(-1,32){\makebox(0,0)[l]{\scriptsize$I^\Diamond$}}

\put(40,54){\makebox(0,0)[t]{\scriptsize$\Box^I$}}
\put(40,66){\makebox(0,0)[b]{\scriptsize$\Diamond^I$}}

\end{picture}
\end{center}
with the two sides omitted being exact replicas of those drawn. In
the horizontal adjunction here, $\Diamond^I$ and $\Box^I$ are just
$\Diamond$ and $\Box$ respectively, and this adjunction is the
original adjunction mentioned at the end of Section~6.

The functors involved in the vertical, contravariant, adjunction,
for $M$ being $\Box$ or $\Diamond$, and $\alpha_A$ a primitive
arrow term of $S5_{\Box\Diamond}$, are defined by
\begin{tabbing}
\hspace{3em}$I^MA=AM$,
\hspace{3em}$I^M\alpha_A=\alpha_{AM}$,\hspace{3em}$I^M(g\cirk
f)=I^Mg\cirk I^Mf$.
\end{tabbing}
That these are indeed functors is guaranteed by the fact that the
equations of $S5_{\Box\Diamond}$ are assumed universally. These
functors will hence exist also when we extend $S5_{\Box\Diamond}$
with new equations, assumed universally. Note that they exist in
the free dyad $S5_{\Box\Diamond}$, but they need not exist in an
arbitrary dyad. (Analogous functors exist in  $S4_\Box$,
$S4_\Diamond$, etc., but they need not exist in arbitrary comonads
and monads.)

\section{Maximality in the context of $S5$}
Consider the following equations, which do not hold in
$S5_{\Box\Diamond}$:
\begin{tabbing}
\hspace{1.7em}\= $\Box\varepsilon^\Box_A=\varepsilon^\Box_{\Box
A}$,\hspace{13em}\= $\varepsilon^\Box_A\cirk\varepsilon^\Box_{\Box
A}\cirk\delta^{\Diamond\Box}_{\Box
A}=\varepsilon^\Box_A\cirk\delta^{\Diamond\Box}_A\cirk\Diamond\varepsilon^\Box_{\Box
A}$,
\\*[1.5ex]
\>
$\Box\varepsilon^\Diamond_A=\delta^{\Box\Diamond}_A\cirk\varepsilon^\Diamond_A\cirk\varepsilon^\Box_A$,\>
$\varepsilon^\Diamond_A\cirk\varepsilon^\Box_A\cirk\delta^{\Diamond\Box}_A=\Diamond\varepsilon^\Box_A$,
\\*[1.5ex]
\> $\Box\varepsilon^\Diamond_{\Diamond
A}\cirk\delta^{\Box\Diamond}_A\cirk\varepsilon^\Diamond_A=\delta^{\Box\Diamond}_{\Diamond
A}\cirk\varepsilon^\Diamond_{\Diamond
A}\cirk\varepsilon^\Diamond_A$,\> $\varepsilon^\Diamond_{\Diamond
A}=\Diamond\varepsilon^\Diamond_A$.
\end{tabbing}
In the left upper corner and the right lower corner we have the
equations ($\varepsilon^\Box$~{\it triv}) and
($\varepsilon^\Diamond$~{\it triv}). The left-hand sides of these
six equations correspond to the six pictures on the left, while
the right-hand sides correspond to the six pictures on the right:
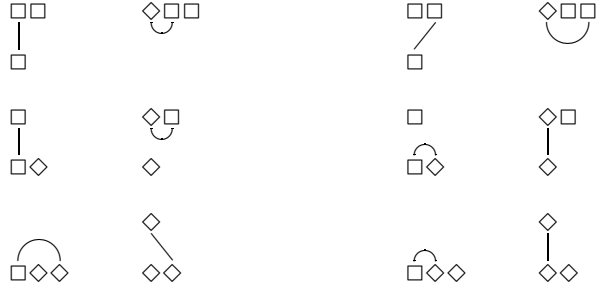
\begin{center}
\begin{picture}(230,110)

\put(0,0){\makebox(0,0)[bl]{$\Box\Diamond\Diamond$}}
\put(11.5,8){\oval(16,16)[t]}

\put(0,40){\makebox(0,0)[bl]{$\Box\Diamond$}}
\put(0,59){\makebox(0,0)[bl]{$\Box$}} \put(4,48){\line(0,1){10}}

\put(0,80){\makebox(0,0)[bl]{$\Box$}}
\put(0,99){\makebox(0,0)[bl]{$\Box\Box$}}
\put(4,88){\line(0,1){10}}

\put(50,0){\makebox(0,0)[bl]{$\Diamond\Diamond$}}
\put(50,19){\makebox(0,0)[bl]{$\Diamond$}}
\put(62,8){\line(-4,5){8}}

\put(50,40){\makebox(0,0)[bl]{$\Diamond$}}
\put(50,59){\makebox(0,0)[bl]{$\Diamond\Box$}}
\put(58,58){\oval(8,8)[b]}

\put(50,99){\makebox(0,0)[bl]{$\Diamond\Box\Box$}}
\put(58,98){\oval(8,8)[b]}

\put(150,0){\makebox(0,0)[bl]{$\Box\Diamond\Diamond$}}
\put(157.5,8){\oval(8,8)[t]}

\put(150,40){\makebox(0,0)[bl]{$\Box\Diamond$}}
\put(150,59){\makebox(0,0)[bl]{$\Box$}}
\put(157.5,48){\oval(8,8)[t]}

\put(150,80){\makebox(0,0)[bl]{$\Box$}}
\put(150,99){\makebox(0,0)[bl]{$\Box\Box$}}
\put(153.5,88){\line(4,5){8}}

\put(200,0){\makebox(0,0)[bl]{$\Diamond\Diamond$}}
\put(200,19){\makebox(0,0)[bl]{$\Diamond$}}
\put(204,8){\line(0,1){10}}

\put(200,40){\makebox(0,0)[bl]{$\Diamond$}}
\put(200,59){\makebox(0,0)[bl]{$\Diamond\Box$}}
\put(204,48){\line(0,1){10}}

\put(200,99){\makebox(0,0)[bl]{$\Diamond\Box\Box$}}
\put(211.5,98){\oval(16,16)[b]}

\end{picture}
\end{center}

The bijections between hom-sets of the horizontal adjunction of
$S5_{\Box\Diamond}$ mentioned at the end of the preceding section
stand behind the horizontal connections in the six pictures on the
left. The same holds when we replace ``horizontal'' by
``vertical'', or ``left'' by ``right''. From that we can conclude
that any of the six equations above when added to
$S5_{\Box\Diamond}$ yields the five remaining ones. Anticipating
matters, we call any of these equations a \emph{preordering
equation} of $S5_{\Box\Diamond}$.

Let $S5_{\Box\Diamond triv}$ be the category defined like
$S5_{\Box\Diamond}$ save that we have as an additional equation
one of the preordering equations of $S5_{\Box\Diamond}$
(universally assumed). To show that $S5_{\Box\Diamond triv}$ is a
preorder, we need to consider first some properties of the functor
$G$ from $S5_{\Box\Diamond}$ to \emph{Gen}.

For every arrow $f$ of $S5_{\Box\Diamond}$, the partition
corresponding to the split equivalence $Gf$ induces a partition on
the occurrences of $\Box$ and $\Diamond$ in the source and target
of $f$, and we call the members of the latter partition the
\emph{equivalence classes} of $f$. An element of an equivalence
class of $f$ is either a \emph{source element} or a \emph{target
element}, and also every such element is either a $\Box$
\emph{element} or a $\Diamond$ \emph{element}.

From the normal form for the arrow terms of $S5_{\Box\Diamond}$ in
the proof of $S5_{\Box\Diamond}$ Coherence in Section~6, we can
conclude that for every arrow $f$ of $S5_{\Box\Diamond}$ the
equivalence classes of $f$ are of one of the following two kinds:
\begin{itemize}
\item[($\Box$)]there is a $\Box$ element that is the rightmost
source element in the class, and is called the \emph{head} of the
class; all the other source elements (if any) are $\Diamond$
elements, and all the target elements (if any) are $\Box$
elements; \item[($\Diamond$)]there is a $\Diamond$ element that is
the rightmost target element in the class, and is called the
\emph{head} of the class; all the other target elements (if any)
are $\Box$ elements, and all the source elements (if any) are
$\Diamond$ elements.
\end{itemize}
Every source $\Box$ element and every target $\Diamond$ element is
a head. Let an element of an equivalence class that is not its
head be called \emph{subordinate}. Every source $\Diamond$ element
and every target $\Box$ element is subordinate. The number of
equivalence classes of an arrow depends only on the type of that
arrow.

Take an arrow ${f\!:A\vdash B}$ of $S5_{\Box\Diamond}$, and
consider an equivalence class $E$ of $f$. For an arbitrary subset
$E'$ of $E$ that contains the head of $E$, there is an arrow
${k_A\!:A'\vdash A}$ built by using essentially
$\varepsilon^\Diamond$, and there is an arrow ${k_B\!:B\vdash B'}$
built by using essentially $\varepsilon^\Box$, such that
${k_B\cirk f\cirk k_A}$ has equivalence classes exactly like $f$
save that $E$ is replaced by $E'$. As a limit case, we may take
$E'$ to be the singleton whose only member is the head of $E$. We
say that ${k_B\cirk f\cirk k_A}$ is obtained by \emph{reducing}
$E$ in $f$ to $E'$. Next we show the following.

\prop{Preorder of $S5_{\Box\Diamond triv}$}{The category
$S5_{\Box\Diamond triv}$ is a preorder, and its skeleton is given
by the following diagram:}

\vspace{-2ex}

\begin{center}
\begin{picture}(60,15)

\put(5,6){\vector(1,0){21}} \put(37,6){\vector(1,0){21}}

\put(0,5){\makebox(0,0){$\Box$}}
\put(64,5){\makebox(0,0){$\Diamond$}}

\put(16,9){\makebox(0,0)[b]{\scriptsize$\varepsilon^\Box$}}
\put(49,9){\makebox(0,0)[b]{\scriptsize$\varepsilon^\Diamond$}}

\end{picture}
\end{center}

\dkz Note first that, for $M_1,M_2\in\{\Box,\Diamond\}$, we have
in $S5_{\Box\Diamond triv}$ that $M_1M_2$ is isomorphic to $M_2$.
To prove these isomorphisms, besides equations we have encountered
previously, we have
\begin{tabbing}
\hspace{9em}$\delta^{\Box M}\cirk\varepsilon^\Box_M$ \=
$=\varepsilon^\Box_{\Box M}\cirk\Box\delta^{\Box M}$,
\hspace{.5em}by ($\varepsilon^\Box$~{\it nat}),
\\*[1.5ex]
\> $=\mj_{\Box M}$, \hspace{.5em}by ($\varepsilon^\Box$~{\it
triv}) and $(\Box M\beta)$,
\end{tabbing}
and we derive analogously
${\varepsilon^\Diamond_M\cirk\delta^{\Diamond M}=\mj_{\Diamond
M}}$. Next, if $M_1,M_2\in\{\Box,\emptyset,\Diamond\}$, then from
$S5_{\Box\Diamond}$ Coherence and the form of the equivalence
classes of the arrows of $S5_{\Box\Diamond}$, we may conclude that
there is at most one arrow from $M_1$ to $M_2$ in
$S5_{\Box\Diamond}$, and hence also in $S5_{\Box\Diamond
triv}$.\qed

\vspace{2ex}

\noindent We define $5S_{\Box\Diamond triv}$ analogously, and
prove in the same manner that it is a preorder, with an isomorphic
skeleton.

The category $S5_{\Box\Diamond}$, as well as $5S_{\Box\Diamond}$,
is maximal in the sense in which $S4_{\Box\Diamond\sharp}$ was
shown maximal in Section~9.

\prop{Maximality of $S5_{\Box\Diamond}$}{The category
$S5_{\Box\Diamond}$ is maximal.}

\dkz Suppose we have the arrow terms ${f_1,f_2\!:A\vdash B}$ of
$S5_{\Box\Diamond}$ such that $f_1=f_2$ does not hold in
$S5_{\Box\Diamond}$. By $S5_{\Box\Diamond}$ Coherence, we have
$Gf_1\neq Gf_2$. Then it can be inferred that there are three
distinct occurrences $x$, $y_1$ and $y_2$ of $\Box$ or $\Diamond$
in $A$ or $B$ such that $y_1$ and $y_2$ are heads of equivalence
classes both in $f_1$ and in $f_2$, and $x$ is in the same class
$E_1$ as $y_1$ in $f_1$, and in the same class $E_2$ as $y_2$ in
$f_2$. So $x$ is a subordinate element both in $f_1$ and $f_2$.

Let $S5^\ast_{\Box\Diamond}$ be obtained by extending
$S5_{\Box\Diamond}$ with $f_1=f_2$, universally assumed. By
reducing $E_1$ in $f_1$ to $\{x,y_1\}$, and every other
equivalence class of $f_1$ to a singleton, we obtain the arrow
$f_1'=k''\cirk f_1\cirk k'$, which is equal in
$S5^\ast_{\Box\Diamond}$ to $f_2'=k''\cirk f_2\cirk k'$. In $f_2'$
the subordinate element $x$ belongs to the equivalence class
$\{x,y_2\}$, while all the other equivalence classes of $f_2'$ are
singletons.

The number $n$ of equivalence classes in $f_1$, $f_2$, $f_1'$ and
$f_2'$ is the same, and we proceed by induction on $n$ to show
that we can derive one of the preordering equations in
$S5^\ast_{\Box\Diamond}$, i.e.\ that $S5^\ast_{\Box\Diamond}$ is
$S5_{\Box\Diamond triv}$. The basis of this induction is when
$n=2$, and then we have cases that are covered by the six
preordering equations of $S5_{\Box\Diamond}$. If $n\geq 4$, then
either $f_1'\cirk C\delta^{\Box M}=f_2'\cirk C\delta^{\Box M}$,
and we can apply the induction hypothesis, or $C\delta^{\Diamond
M}\cirk f_1'=C\delta^{\Diamond M}\cirk f_2'$, and we can apply the
induction hypothesis (here $C$ is a modality, possibly empty).

If $n=3$, we proceed either as when $n\geq 4$, or we have an
additional case in which we rely also on the vertical adjunction
involving $I^\Box$ and $I^\Diamond$, which obtains also in
$S5^\ast_{\Box\Diamond}$ (see the preceding section). For example,
if we find ourselves in the situation that corresponds to the
following pictures:
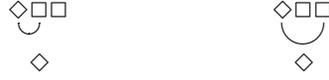
\begin{center}
\begin{picture}(120,30)

\put(8,0){\makebox(0,0)[bl]{$\Diamond$}}
\put(0,20){\makebox(0,0)[bl]{$\Diamond\Box\Box$}}
\put(8,18){\oval(8,8)[b]}

\put(108,0){\makebox(0,0)[bl]{$\Diamond$}}
\put(100,20){\makebox(0,0)[bl]{$\Diamond\Box\Box$}}
\put(111.5,18){\oval(16,16)[b]}

\end{picture}
\end{center}
by the vertical adjunction, we pass first to
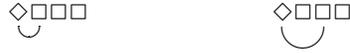
\begin{center}
\begin{picture}(120,20)

\put(0,10){\makebox(0,0)[bl]{$\Diamond\Box\Box\Box$}}
\put(8,8){\oval(8,8)[b]}

\put(100,10){\makebox(0,0)[bl]{$\Diamond\Box\Box\Box$}}
\put(111.5,8){\oval(16,16)[b]}

\end{picture}
\end{center}
and then by precomposing with $\Diamond\Box\delta^{\Box\Box}$ we
obtain
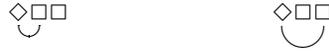
\begin{center}
\begin{picture}(120,20)

\put(0,10){\makebox(0,0)[bl]{$\Diamond\Box\Box$}}
\put(8,8){\oval(8,8)[b]}

\put(100,10){\makebox(0,0)[bl]{$\Diamond\Box\Box$}}
\put(111.5,8){\oval(16,16)[b]}

\end{picture}
\end{center}
i.e.\ the preordering equation in the right upper corner. This is
enough to show that $S5_{\Box\Diamond}$ is maximal.\qed

\vspace{2ex}

\noindent The category $5S_{\Box\Diamond}$ is shown to be maximal
in the same manner.

It is shown in \cite{D99} (\emph{Addenda and Corrigenda}, Section
5.11) that the maximality of comonads, i.e.\ of $S4_\Box$, entails
an analogous maximality of adjunction. In the same way, the
maximality of $S5_{\Box\Diamond}$ or $5S_{\Box\Diamond}$ entails
the maximality of trijunction, as we will show below. We cannot
extend this notion with new equations in the canonical language of
trijunctions, equations being assumed universally (cf.\
Section~9), without trivializing the notion: any equation in the
canonical language will hold.

To infer the maximality of adjunction from the maximality of
comonads, or the maximality of monads, we can proceed not as in
the reference mentioned above, but by appealing to the square of
adjunctions of the preceding section. The category ${\cal A}^{\cal
A}$ corresponds to the comonad, and ${\cal B}^{\cal B}$ to the
monad. Any arrow of the freely generated adjunction is in one of
four disjoint categories, which correspond to the categories
$C{\cal A}^{\cal B}$, $C{\cal B}^{\cal B}$, $C{\cal A}^{\cal A}$
and $C{\cal B}^{\cal A}$ (see the preceding section). By the
horizontal and vertical adjunctions, any such equation can be
reduced to a new equation of comonads or monads.

There is a square of trijunctions analogous to the square of
adjunctions. Suppose we have a trijunction given by the categories
$\cal A$ and $\cal B$, a functor $U$ from $\cal A$ to $\cal B$,
and the functors $L$ and $R$ from $\cal B$ to $\cal A$, with $L$
being left adjoint and $R$ right adjoint to $U$. Then, with arrows
of right adjoints being more finely dotted, we have
\begin{center}
\begin{picture}(80,100)(0,-20)

\put(-11,13){\vector(0,1){40}}
\multiput(-4,53)(0,-4){4}{\line(0,-1){2}}
\multiput(-4,26)(0,-4){3}{\line(0,-1){2}}
\put(-4,15){\vector(0,-1){2}}
\multiput(3,13)(0,2){20}{\circle*{.1}} \put(3,51){\vector(0,1){2}}

\put(74,13){\vector(0,1){40}}
\multiput(81,53)(0,-4){4}{\line(0,-1){2}}
\multiput(81,26)(0,-4){3}{\line(0,-1){2}}
\put(81,15){\vector(0,-1){2}}
\multiput(88,13)(0,2){20}{\circle*{.1}}
\put(88,51){\vector(0,1){2}}

\put(69,70){\vector(-1,0){58}}
\multiput(11,63)(4,0){6}{\line(1,0){2}}
\multiput(46,63)(4,0){5}{\line(1,0){2}}
\put(67,63){\vector(1,0){2}}
\multiput(69,56)(-2,0){28}{\circle*{.1}}
\put(13,56){\vector(-1,0){2}}

\put(69,8){\vector(-1,0){58}}
\multiput(11,1)(4,0){6}{\line(1,0){2}}
\multiput(46,1)(4,0){5}{\line(1,0){2}} \put(67,1){\vector(1,0){2}}
\multiput(69,-6)(-2,0){28}{\circle*{.1}}
\put(13,-6){\vector(-1,0){2}}

\put(-1,2){\makebox(0,0){${\cal A}^{\cal A}$}}
\put(82,2){\makebox(0,0){${\cal B}^{\cal A}$}}
\put(-1,62){\makebox(0,0){${\cal A}^{\cal B}$}}
\put(82,62){\makebox(0,0){${\cal B}^{\cal B}$}}

\put(-11,32){\makebox(0,0)[r]{\scriptsize${\cal A}^R$}}
\put(-3,32){\makebox(0,0){\scriptsize${\cal A}^U$}}
\put(4,32){\makebox(0,0)[l]{\scriptsize${\cal A}^L$}}

\put(74,32){\makebox(0,0)[r]{\scriptsize${\cal B}^R$}}
\put(82,32){\makebox(0,0){\scriptsize${\cal B}^U$}}
\put(89,32){\makebox(0,0)[l]{\scriptsize${\cal B}^L$}}

\put(40,-8){\makebox(0,0)[t]{\scriptsize$R^{\cal A}$}}
\put(40,2){\makebox(0,0){\scriptsize$U^{\cal A}$}}
\put(40,11){\makebox(0,0)[b]{\scriptsize$L^{\cal A}$}}

\put(40,54){\makebox(0,0)[t]{\scriptsize$R^{\cal B}$}}
\put(40,64){\makebox(0,0){\scriptsize$U^{\cal B}$}}
\put(40,73){\makebox(0,0)[b]{\scriptsize$L^{\cal B}$}}

\end{picture}
\end{center}
The category ${\cal B}^{\cal B}$ here corresponds to dyads, i.e.\
$S5_{\Box\Diamond}$, and ${\cal A}^{\cal A}$ to codyads, i.e.\
$5S_{\Box\Diamond}$. Any arrow of the freely generated trijunction
is in one of four disjoint categories, which correspond to the
four categories in the square of trijunctions above. For example,
to ${\cal A}^{\cal B}$ there corresponds a category $C{\cal
A}^{\cal B}$ whose objects are $L$, $R$, $LUL$, $RUL$, etc., to
${\cal B}^{\cal B}$ there corresponds a category $C{\cal B}^{\cal
B}$ whose objects are $\emptyset$, $UL$, $UR$, $ULUL$, $URUL$,
etc., to ${\cal A}^{\cal A}$ there corresponds a category $C{\cal
A}^{\cal A}$ whose objects are $\emptyset$, $LU$, $RU$, $LULU$,
$LURU$, etc., and, finally, to ${\cal B}^{\cal A}$ there
corresponds a category $C{\cal B}^{\cal A}$ whose objects are $U$,
$ULU$, $URU$, $ULULU$, $ULURU$, etc. Here, $\emptyset$ corresponds
to identity functors.

By these horizontal and vertical adjunctions, any new equation of
trijunctions can be reduced to a new equation of dyads or codyads.
So the maximality of trijunction can be inferred from the
maximality of $S5_{\Box\Diamond}$, or the maximality of
$5S_{\Box\Diamond}$.

To make this inference, we could also proceed as in \cite{D99}
(\emph{Addenda and Corrigenda}, Section 5.11). The category
$C{\cal A}^{\cal B}$ is isomorphic by the functor $U$ to a
subcategory ${\cal B}'$ of $C{\cal B}^{\cal B}$, and $C{\cal
B}^{\cal B}$ and ${\cal B}'$ together with the functors $LU$, $RU$
and the inclusion functor from ${\cal B}'$ to $C{\cal B}^{\cal B}$
make a trijunction isomorphic to the original trijunction. The
category ${\cal B}'$ is isomorphic to the category $(C{\cal
B}^{\cal B})^{UL}_{UR}$ of Section~8. Any new equation for
trijunctions corresponds by this isomorphism to a new equation of
dyads.

We will not consider here the extension of $S5_{\Box\Diamond}$
with the arrows $\chi_A^{\Diamond\Box}$ or $\chi_A^{\Box\Diamond}$
of Section~5. With $\chi_A^{\Box\Diamond}$ we would obtain a
$\Box\Diamond$-structure that is both $S5_{\Box\Diamond}$ and
$5S_{\Box\Diamond}$, at the same time. With this structure, we
come close to the Frobenius monads of \cite{LAW69} (pp.\ 151-152);
namely, dyads where $\Box$ and $\Diamond$ coincide, and where
$\delta^{\Box\Diamond}$ and $\delta^{\Diamond\Box}$ coincide
respectively with $\delta^{\Box\Box}$ and
$\delta^{\Diamond\Diamond}$ (alternatively, these are codyads
where $\Box$ and $\Diamond$ coincide). We deal with them in
\cite{DP08a}.

\vspace{1.8ex}

\noindent {\small {\it Acknowledgements$\,$}. Work on this paper
was supported by the Ministry of Science of Serbia. We are very
grateful to an anonymous referee for his good will and effort in
reading our paper, and for making a useful suggestion concerning
presentation.}

\end{document}